\documentclass{article}
\usepackage[utf8]{inputenc}
\usepackage{graphicx}
\usepackage[hyperref,dvipsnames]{xcolor}

\usepackage{amsmath, color, verbatim}
\usepackage[colorlinks=true]{hyperref}

\usepackage{amsfonts, mathtools, enumitem}
\usepackage{todonotes}
\usepackage{titlesec}
\usepackage{algorithm}
\usepackage{algorithmic}
\usepackage{caption}
\usepackage{placeins}
\hypersetup{
  colorlinks=true,
  citecolor=MidnightBlue,
  linkcolor=purple!80
  }
\usepackage{amssymb}
\usepackage{amsthm}
\usepackage{bbm}
\usepackage[left=1in,right=1in,top=1in,bottom=1in]{geometry}
\usepackage{svg}
\usepackage{thm-restate}

\makeatletter
\@ifundefined{newcounteralias}{}{%
  \renewcommand\thmt@autorefsetup{%
    \@xa\def\csname\thmt@envname autorefname\@xa\endcsname\@xa{\thmt@thmname}%
  }%
}
\makeatother

\usepackage[nameinlink,noabbrev,capitalize]{cleveref}

\DeclareMathOperator*{\esssup}{ess\,sup}

\newcommand{\norm}[2]{\left\lVert #1\right\rVert_{#2}}

\newtheorem{theorem}{Theorem}[section]

\newtheorem{lemma}[theorem]{Lemma}
\newtheorem{prop}[theorem]{Proposition}

\newtheorem{ass}[theorem]{Assumption}

\newtheorem{remark}[theorem]{Remark}

\crefname{ass}{Assumption}{Assumptions}
\crefname{prop}{Proposition}{Propositions}

\crefname{remark}{Remark}{Remarks}

\newcommand{\grad}{\nabla}

\newcommand{\weaklyto}{\rightharpoonup}
\newcommand{\weakstar}{\stackrel{*}\rightharpoonup}
\newcommand{\cts}{\hookrightarrow}

\newcommand{\ctsCompact}{\xhookrightarrow{c}}
\newcommand{\E}{\mathbb{E}}
\newcommand{\pP}{\mathbb{P}}
\newcommand{\R}{\mathbb{R}}

\newcommand{\cR}{\mathcal{R}}
\newcommand{\cJ}{\mathcal{J}}

\newcommand{\CVaR}{\textup{CVaR}}
\newcommand{\D}{\textup{ d}}
 \makeatletter
\patchcmd{\l@section}{1.0em}{0.2em}{}{}
\makeatother

\begin{document}

\title{Optimality conditions for risk-averse control of elliptic variational inequalities under uncertainty
}

\author{Amal Alphonse\thanks{Faculty of Mathematics,
  University of Duisburg-Essen, Thea-Leymann-Straße 9, 
  45127 Essen, Germany. Email: \url{amal.alphonse@uni-due.de}}}
\date{}
\newcommand{\skipeqs}[1]{\addtocounter{equation}{#1}}

\maketitle
\begin{abstract}
We establish first-order optimality conditions for stochastic optimisation problems governed by elliptic obstacle variational inequalities with random data and convex risk measures. The nonsmoothness of both the control-to-state mapping and the risk functional presents substantial difficulties in deriving stationarity conditions. We derive a full C-stationarity system on an almost-sure (i.e., scenario-wise) level for local minimisers under either monotonicity of the risk measure (without requiring any differentiability), or under Fréchet differentiability of the risk measure at the optimal loss. Our results considerably improve previously available stationarity conditions by establishing almost-sure sign and contact orthogonality conditions, together with characterisations of the support of the adjoint multipliers, leading to a complete C-stationarity optimality system. A key contribution is the construction of an unweighted adjoint that is compatible with the risk-weighted optimality system, even in scenarios assigned zero risk weight. Our theory applies to a wide range of risk measures, including both monotone and non-monotone as well as smooth and nonsmooth ones. We illustrate the results for a mean--CVaR optimal control problem.%
\end{abstract}

 \tableofcontents

\section{Introduction}\label{sec:problem-setting}
In this paper, we study first-order optimality conditions for the stochastic optimisation problem
\begin{equation}\label{eq:OCProblem}
\min_{u\in U_{ad}} \cR[\cJ(S(u))]+\varrho(u),
\end{equation}
where $U_{ad}$ is an admissible set of controls, $\cJ$ is an objective function, $\cR$ is a risk measure, $\varrho$ is the control cost, and $S$ is the solution map of a random elliptic variational inequality (VI). Indeed, for almost every $\omega\in\Omega$ taken from a probability space $(\Omega, \mathcal{F}, \mathbb{P})$, the object $S(u)(\omega)=:y(\omega)$ is the unique solution of the obstacle problem%
\begin{equation}\label{eq:VIProblemFormal}
\begin{aligned}
&y(\omega)\in V,\qquad y(\omega)\leq\psi(\omega), \qquad \langle A(\omega)y(\omega)-f(\omega)-B(\omega)u,\,y(\omega)-v\rangle\leq0
\qquad\forall v\in V:\ v\leq\psi(\omega),
\end{aligned}
\end{equation}
where $V$ is either $H^1(D)$ or $H^1_0(D)$ on a bounded domain $D \subset \mathbb{R}^d$. %
Precise assumptions on these objects and spaces will be introduced in due course.

The aim of this work is to derive a C-stationarity system for this problem, addressing a gap in the literature. This paper is closely related to \cite{AGHS}, where the very same problem  is considered and various forms of ($\mathcal{E}$-almost) \textit{weak} stationarity systems are derived. 
In particular, %
the stationarity result in \cite{AGHS} required a
single-valued subdifferential and (continuous) differentiability 
of $\cR$ on all of $L^p(\Omega)$, together with an
additional restriction on the integrability exponents.
Even under these assumptions, the resulting optimality
systems retained weak or $\mathcal E$-almost conditions,
or approximating (regularisation-dependent) statements necessitated by technical difficulties, and thus these systems fall short of what is known in (and what we expect from) the deterministic setting. In this work, we will overcome all of these issues under much weaker assumptions on the risk measure and (to the best of our knowledge) present the first stationarity system that corresponds fully to well known and established notions of stationarity in PDE- and VI-constrained and nonsmooth optimisation. Indeed, we will obtain C-stationarity under the assumption that \textit{either} $\cR$ is a monotone risk measure (thus allowing for nonsmooth risk measures), \textit{or} that $\cR$ is Fr\'echet differentiable at the optimal loss $\cJ(y^*)$ (thus allowing for potentially non-monotone risk functionals) without requiring differentiability
everywhere. %

Due to this generality, our results cover a wide range of risk measures in practice, including all finite-valued coherent and all convex monetary risk measures. In particular, our results apply to the widely-used conditional value-at-risk (CVaR) and mean--CVaR risk 
measures, which are monotone but generally nonsmooth. For these important examples, we obtain a complete C-stationarity system.

The key developments that enable this strengthening of the results are the following:
\begin{itemize}[itemsep=0.1pt]
    \item the tools developed in the recent works \cite{AHRW, AW} enable us to prove strong a.s. orthogonality properties as well as avoid requiring Egorov-type approximations, removing the need for the so-called $\mathcal{E}$-almost qualifiers that were needed in \cite{AGHS};
    \item a Minty-type argument combined with Mazur's lemma enables us to establish almost-sure sign conditions despite having only the weak convergence of the penalised adjoints and risk multipliers,  %
    yielding  in particular complementarity-type behaviour of monotone risk measures and an adjoint statement, completing the C-stationarity conditions, \textit{without} requiring any differentiability on $\cR$.
\end{itemize}
Unlike weak stationarity, C-stationarity incorporates
additional sign and complementarity conditions on the
adjoint variables and multipliers, providing a sharper
first-order description of local optimality. These
conditions are also relevant from a numerical perspective,
as they characterise the optimality structure expected
in the limits of penalisation and regularisation schemes
and provide additional criteria for assessing numerical
approximations.

\paragraph{Contributions} 
As explained, the main contribution of this work is the derivation of a complete almost-sure C-stationarity system for arbitrary local minimisers of the risk-averse obstacle control problem. Our optimality results hold for two different classes of risk measures $\cR$:
\begin{itemize}[itemsep=0.1pt]
    \item risk functionals that are Fr\'echet differentiable
    at the optimal loss $\cJ(y^*)$;
    \item monotone risk functionals (which can be nonsmooth).
\end{itemize}
Compared with \cite{AGHS}, we strengthen the available
optimality conditions by establishing almost-sure sign
conditions, separate contact orthogonality for the positive
and negative parts of the adjoints, and exact inactive-set
support conditions without approximation by exceptional
sets. A central difficulty is that the weak convergence of
risk multipliers and penalised adjoints does not generally
preserve their product structure. We overcome this
difficulty by constructing an unweighted adjoint compatible
with the risk-weighted optimality system, including in
scenarios where the limiting risk multiplier vanishes.
Crucially, this construction requires no differentiability
of the risk functional when it is monotone. Alternatively,
without monotonicity, Fr\'echet differentiability at the
optimal loss suffices. 
We also remove an additional exponent restriction and a weak lower semicontinuity
assumption on the control-cost derivative imposed
in \cite{AGHS}.
\paragraph{Related literature}
Elliptic variational inequalities of obstacle type arise in numerous applications in physics, biology and economics such as 
contact mechanics, biological processes in cells, ecology, fluid flow, finance, and other areas where free boundary problems appear. Uncertainty in the coefficients, source terms and obstacles is more realistic  in many applications where one does not have access to exact or precise data or measurements, but knows their general distribution. Their analysis (under uncertain data terms) has been investigated in
\cite{Gwinner2000,Gwinner2006a,Forster2010,
Kornhuber2014,Bierig2015}.
A related body of work concerns stochastic variational
inequalities and stochastic mathematical programs with
equilibrium constraints; see, for example,
\cite{Patriksson1999,Evgrafov2004a,
shapiro2006stochastic,Shapiro2008,
Rockafellar2016,Rockafellar2018,Rockafellar2020}.
These approaches encompass different treatments of
uncertainty and equilibrium constraints. In the present
problem, as in \cite{AGHS}, the control is a deterministic here-and-now
decision, while the state is a random function satisfying
the elliptic variational inequality almost surely.

For deterministic optimal control problems governed
by elliptic obstacle variational inequalities, an extensive
literature has developed on sensitivity analysis,
generalised derivatives, optimality conditions and
numerical methods; see for example
\cite{Barbu1984,MR2515801,MR2822818,Wachsmuth,MR3796767,
RW}.
In particular, the nonsmoothness of the obstacle
solution operator motivates different stationarity
concepts, including weak, C-, M- and strong
stationarity. More recently, \cite{AW}
studied penalty approximations of the obstacle problem
and established convergence properties of the
derivatives of the penalised solution operators.
For elliptic obstacle-type quasi-variational
inequalities, \cite{AHRW} established differentiability properties of the minimal and maximal solution maps and derived C-stationarity conditions for associated deterministic optimal control problems. 

Risk measures provide a flexible framework for incorporating the distribution of uncertain outcomes into optimisation criteria; see \cite{Shapiro2009,Pflug2007}. Their use in PDE-constrained optimisation under uncertainty has been developed in, among others, \cite{Kouri2018,Kouri2019,Kouri2019a}. Recent contributions include for example the treatment of entropic risk measures in parabolic PDE-constrained control \cite{GuthEtAl} and stochastic optimal design under risk measures \cite{AKV}. Related developments in risk-averse stochastic optimisation and stochastic variational inequalities, particularly concerning sample-based approximations, can be found in e.g. \cite{MelnikovMilz}. In contrast to the present problem, these works do not address the simultaneous nonsmoothness arising from an elliptic obstacle solution operator and a potentially nonsmooth risk functional.

Closer to the present setting,
\cite{hertlein2022inexact} investigated an
inexact bundle method for optimal control problems governed by
deterministic and stochastic obstacle problems with the risk measure chosen to be the expectation.
The risk-averse problem \eqref{eq:OCProblem} with random elliptic
variational inequality constraints was subsequently
studied in \cite{AGHS}, as already mentioned, where penalisation
was used to derive weak forms of stationarity conditions. The present work is inspired by \cite{AGHS} and builds upon the foundation laid therein.%

\paragraph{Outline of paper} In \cref{sec:main_stuff}, we start  setting the scene by introducing necessary assumptions and problem data in order for \eqref{eq:OCProblem} to be well defined. We state and discuss our main result \cref{thmmain} on the first-order C-stationarity conditions, compare it with the literature, and provide some examples satisfying our framework.  We also illustrate the optimality system for a mean--CVaR control problem, deriving an explicit characterisation of the risk multipliers and a projection formula for the optimal control. 
In \cref{sec:optimisation-problem} we start the preparation for proving our main theorem by formulating the regularisation/penalisation scheme, where the variational inequality \eqref{eq:VIProblemFormal} (which appears via $S$ in \eqref{eq:OCProblem}) is approximated by semilinear partial differential equations, and some initial convergence results are given; much of this is taken from \cite{AGHS}, but we weaken certain assumptions from the cited paper. Next in \cref{sec:strengthened}, we prove new orthogonality conditions for the multipliers, making use of some methodology in \cite{AHRW, AW} and a careful limiting analysis, culminating in an enhanced weak stationarity system \cref{thmenhancedWeakStationarity}. In addition, we lay the groundwork to prove C-stationarity by developing the technical tool in \cref{lem:technical_tool}. In \cref{sec:C_stat_Frechet} we prove the C-stationarity result \cref{thmmain} under Fréchet differentiability of $\cR$, and in \cref{sec:factorisation} we prove it under monotonicity of $\cR$. Finally, we finish the paper with some concluding comments in \cref{sec:conclusion}.

\section{Problem formulation, main results, and examples}\label{sec:main_stuff}
In this section, we begin by formulating the quantities and spaces in the problem of interest \eqref{eq:OCProblem}, which we recall here is
\[\min_{u\in U_{ad}} \cR[\cJ(S(u))]+\varrho(u),\]
and stating the standing assumptions needed for this paper. We will then discuss our main result before moving onto some examples and interpreting  our optimality conditions for a mean--CVaR application. 
\subsection{Notation and standing assumptions}\label{sec:assumptions}
Let us proceed by first clarifying some notation.
\paragraph{Notation}%
Given a Banach space $X$ and a probability space $(\Omega,\mathcal{F}, \pP)$, for exponents $t \geq 1$, the Bochner space $L^t(\Omega; X):=L^t(\Omega,\mathcal{F}, \pP; X)$ is the set of all (equivalence classes of) strongly measurable functions $y\colon \Omega \rightarrow X$ %
such that the following norm is finite:
\begin{equation*}
\lVert y \rVert_{L^t(\Omega;X)}:= \begin{cases}
                                     (\int_\Omega \lVert y(\omega) \rVert_X^t \D \pP(\omega))^{1/t}  &\text{ for } t \in [1, \infty),\\
                                    \esssup_{\omega \in \Omega} \lVert y(\omega) \rVert_X &\text{ for } t=\infty.
                                    \end{cases}
\end{equation*}
We set $L^t(\Omega)=L^t(\Omega;\R)$ for the space of random variables with finite $t$-moments.  The expectation of  a random variable $Z\colon \Omega \rightarrow \R$ is defined by $\E[Z]:=\int_{\Omega} Z(\omega)\D \pP(\omega).$
Given Banach spaces $X,Y$, the notation $\mathcal{L}(X,Y)$ denotes the set of bounded and linear maps from $X$ into $Y$. A map $A\colon \Omega \to \mathcal{L}(X,Y)$ is said to be \textit{uniformly measurable} in $\mathcal{F}$ if there exists a sequence of countably-valued operator random variables in $\mathcal{L}(X,Y)$ converging almost everywhere to $A$ in the uniform operator topology. 

Recall that the \textit{subdifferential} of a convex function $F\colon X \to \mathbb R$ at $z$ is the set $\partial F(z) \subset X^*$ defined as
\[\partial F(z) := \{ g \in X^* : F(z)-F(x) \leq \langle g, z-x\rangle_{X^*,X} \quad \forall x \in X\}.\]
Given a constant $t \in [1,\infty)$, we use $t'$ to denote its H\"older conjugate, so that ${1}\slash {t}+ {1}\slash{t'} = 1$. We write $\cts$ to mean a continuous embedding and $\ctsCompact$ for a compact embedding. Whenever we write the duality pairing $\langle \cdot, \cdot \rangle$ without specifying the spaces,  we mean $\langle \cdot, \cdot \rangle_{V^*,V}$.

\paragraph{Standing assumptions}Let us now describe the basic problem setup more precisely:
\begin{enumerate}[label=(\roman*)]\itemsep0em 
    \item $D \subset \mathbb{R}^d$ is a bounded Lipschitz domain for $d \leq 4$,  and take 
    \[H:=L^2(D) \quad\text{ and } V \in \{H^1(D), H^1_0(D)\}.\] 
    \item $U_{ad} \subset U$ is a non-empty, closed and convex set where the control space $U$ is a Hilbert space. 
    \item $(\Omega,\mathcal{F}, \pP)$ is a complete probability space.%
    \item $f \in L^r(\Omega;V^*)$ is a source term and $\psi \in L^s(\Omega;V)$ is an obstacle  for $r,s \in [2,\infty]$.
\item  
{$\cR\colon L^p(\Omega) \to \R$ for $p \in [1,\infty)$ and $\varrho\colon U \to \R$.} 
\end{enumerate}
We make the following assumption on  the operators in \eqref{eq:VIProblemFormal}. 
\begin{ass}[On the operators]\label{ass:newCombinedOnAandB}
The operators $A\colon \Omega \rightarrow \mathcal{L}(V,V^*)$ and $B\colon \Omega \rightarrow \mathcal{L}(U,V^*)$ are uniformly measurable and there exist positive constants $C_a, C_b, C_c$ such that for all $y, z \in V$ and $u\in U$ and a.s.~in $\omega \in \Omega$,
\begin{equation}\label{ass:boundedAndCoercive}
\begin{aligned}
     \langle A(\omega)y,y\rangle   &\geq C_a \lVert y\rVert_V^2,\\   
    \langle A(\omega)y, z \rangle   &\leq C_b\norm{y}{V}\norm{z}{V},\\
     \langle B(\omega)u, z \rangle  &\leq C_c\norm{u}{U}\norm{z}{V}.
\end{aligned} 
\end{equation}
Furthermore we assume
\begin{equation}
\text{If $u_n \weaklyto u$ in $U$ then $B(\omega)u_n \to B(\omega)u$ in $V^*$ a.s.}
\label{ass:BcompletelyCts}\end{equation}
\end{ass}
In examples it is typical that the embedding $U \ctsCompact V^*$ is compact and the complete continuity assumption on $B$ above is designed to mirror this compact embedding. %
As noted in \cite{AGHS}, the operators \(A\) and \(B\) may be generated by random fields. The uniform boundedness and coercivity assumptions are satisfied, for example, by suitable parameterisations involving bounded random variables.

Regarding the objective function $\cJ$ in \eqref{eq:OCProblem}, we assume that it is the superposition of a function $J\colon V \times \Omega \to \R$, i.e., for $y\colon \Omega \rightarrow V,$ we define the superposition operator $\mathcal{J}(y)\colon \Omega \rightarrow \R$ by 
\[\cJ(y)(\omega) := J(y(\omega),\omega).\]
We make the following assumptions on the objective integrand $J$. These ensure properties such as integrability and  Fr\'echet differentiability.

\begin{ass}[On the objective $J$]\label{ass:onJ}
Let $q$ be an exponent satisfying
 \begin{equation}
2 \leq q < \infty, \qquad p < q \leq \min(r,s),  \label{eq:defnOfq}
\end{equation}
and define 
\[\tilde p:= \frac{pq}{q-p}.\]
Assume the following:
\begin{enumerate}[label=(\roman*)]
\item $J\colon V \times \Omega \rightarrow \R$ is a Carath\'eodory  function\footnote{That is, $J(v,\cdot)$ is measurable for {every} $v$ and $J(\cdot,\omega)$ is continuous for {$\pP$-a.e.}~$\omega.$}  and  there exists $C_1 \in L^p(\Omega)$ and $C_2 \geq 0$ such that
\begin{equation*}
    |J(v,\omega)| \leq C_1(\omega) + C_2\norm{v}{V}^{q\slash p}.
\end{equation*}
\item\label{item:JandJyCC}  $J(\cdot,\omega)\colon V \to \mathbb{R}$ is continuously Fr\'echet differentiable for a.e. $\omega \in \Omega$ and $J_y(\cdot, \cdot)\colon V \times \Omega \to V^*$ is Carath\'eodory. 
\item\label{ass:finiteCase}
There exists $\tilde{C}_1 \in L^{\tilde{p}}(\Omega)$ and $\tilde{C}_2 \geq 0$ such that
\begin{equation}
\norm{J_y(v,\omega)}{V^*} \leq \tilde{C}_1(\omega) + \tilde{C}_2\norm{v}{V}^{q/\tilde{p}}.
\label{ass:JyBoundedFinite}
\end{equation}
\end{enumerate}
\end{ass}
Regarding the exponents introduced above, note the relationship
\begin{equation}\label{eq:qprimeRelatedToPprimeEtc}
\frac1{q'}=\frac1{p'}+\frac1{\tilde p}, \qquad \tilde p > p, \qquad \tilde p \geq q'.
\end{equation}
Thanks to \cref{ass:onJ}, the superposition  operator $\cJ\colon L^q(\Omega;V) \to L^p(\Omega)$ is continuously Fr\'echet differentiable with $\mathcal{J}'(y)(h)
= J_y(y(\cdot), \cdot)h(\cdot)$, and $\mathcal{J}'\colon L^q(\Omega;V) \to L^{\tilde p}(\Omega;V^*)$ is continuous, see \cite[\S 3.2]{AGHS}.

Finally, we make the following assumptions on the risk measure and the regularisation appearing in \eqref{eq:OCProblem}.
\begin{ass}[On $\cR$ and $\varrho$]\label{ass:onvarrho}
Assume the following:
\begin{enumerate}[label=(\roman*)]\itemsep0em 
\item  {$\cR\colon L^p(\Omega) \to \R$} is convex and lower semicontinuous.%
\item\label{item:varrhoConvexGateaux} $\varrho\colon U \to \mathbb{R}$ is convex and G\^{a}teaux  differentiable. %

\end{enumerate}
\end{ass}

Due to \cref{ass:onvarrho}, it follows that $\cR$ is continuous \cite[Proposition 2.111]{Bonnans2013} and Hadamard directionally differentiable \cite[Proposition 2.126 (v)]{Bonnans2013} with
\[\cR'[z](h) = \sup_{\nu \in \partial\cR[z]} \E[\nu h ].\]
Note also that \cref{ass:onvarrho} \ref{item:varrhoConvexGateaux} implies weak lower semicontinuity of $\varrho$.

\subsection{Main results and discussion}
Our main result, on C-stationarity optimality for \eqref{eq:OCProblem}, requires (in addition to the standard standing assumptions above) that $\cR$ is either Fréchet differentiable at the optimal loss $\cJ(y^*)$, or that $\cR$ is monotone, where by monotonicity we mean $X\leq Y$ a.s. implies $\cR[X]\leq\cR[Y]$. In the first scenario, this is a strictly \textit{less} restrictive assumption than the comparable result in \cite[Corollary 4.12]{AGHS} and our theory proves even stronger statements.  The ability to treat the second case of a nonsmooth, monotone risk measure, is the main novelty of this work and requires new tools. 

Before stating the theorem, we briefly recall some notions related to capacity (and refer the reader to \cite[\S 6.4.3]{MR1756264} and \cite[Section~3]{MR3796767} for full details). Since we work with $V \in \{H^1_0(D), H^1(D)\}$, $V$ is a Dirichlet space. This Dirichlet space structure enables one to define notions of capacity and quasi-continuity (heuristically, capacity is a way to measure sets finer than through the Lebesgue measure), see \cite[\S 2.1]{FukushimaBook} or \cite[Section 2]{Wachsmuth2014:2} for the $H^1_0(D)$ setting. 

Since $\psi(\omega) - y^*(\omega) \in V$ a.s., it has a quasi-continuous representative and we will identify it with its representative. Defining the set
\[
\widehat D :=
\begin{cases}
D, & \text{if } V=H^1_0(D),\\
\overline D, & \text{if } V=H^1(D),
\end{cases}
\]
we can define for a.e. $\omega$ the \textit{coincidence} or \textit{active} set
\begin{equation}
\mathcal A^*(\omega):=\{x\in \widehat{D}:y^*(\omega,x)=\psi(\omega,x)\}.\label{eq:coincidence_set}    
\end{equation}
This set is defined up to sets of capacity zero. Finally, we say that a statement holds \textit{quasi-everywhere} and write q.e. for short if it holds everywhere except on a set of capacity zero.

We are ready to state the main result. To simplify notation, we suppress the explicit dependence
on $\omega$ in the following optimality conditions (and all identities and inequalities involving random quantities
are understood to hold almost surely unless stated otherwise).
\begin{theorem}[C-stationarity]%
\label{thmmain}
Let  \cref{ass:newCombinedOnAandB,ass:onJ,ass:onvarrho} hold, and let $u^*\in U_{ad}$ be a local minimiser of \eqref{eq:OCProblem} with corresponding state $y^*=S(u^*)$. In addition, assume that either $\cR$ is monotone or that $\cR$ is Fréchet differentiable at $\cJ(y^*)$.

Then there exist
\[
\zeta^*\in L^q(\Omega;V^*),\qquad
\pi^*\in L^{p'}(\Omega),\qquad
q^*\in L^{q'}(\Omega;V),\qquad
\lambda^*\in L^{q'}(\Omega;V^*),
\]
such that
\begin{subequations}\label{eq:factorisedSystem}
\begin{align}
Ay^*-f-Bu^*+\zeta^*&=0  &&\text{a.s.},\label{eq:stateCNew}\\
\zeta^*\geq0,
\quad y^*\leq\psi,
\quad\langle\zeta^*,y^*-\psi\rangle&=0 &&\text{a.s.},\label{eq:stateCompCNew}\\
\pi^* &\in\partial\cR[\cJ(y^*)], \label{eq:riskCNew}\\
        A^*q^*+\lambda^*
-\pi^*J_y(y^*,\cdot)&=0
&&\text{a.s.},\label{eq:system_qstar_eqn}\\
\E[\langle B^*q^*,u^*-v\rangle_{U^*,U}]
+\varrho'(u^*)(u^*-v)
&\le0 &&\forall v\in U_{ad},\label{eq:system_q_star_control_ineq}
\\
\langle \lambda^*, q^*\rangle &\geq 0 &&\text{a.s.},\label{eq:lambdaqgeq0Nonsmooth}\\
\langle\zeta^*,(q^*)^+\rangle
=\langle\zeta^*,(q^*)^-\rangle
&=0 &&\text{a.s.},\label{eq:zeta_q_pos_neg_zero}
\\
\langle\lambda^*,v\rangle
&=0
&& 
\forall v\in V:\ v=0\ \text{q.e. on }\mathcal A^*, 
\text{ a.s.},\label{eq:lambdastar_orth_qw}
\end{align}
\end{subequations}
with $\pi^* \geq 0$ a.s. if $\cR$ is monotone and $\pi^*=\mathcal R'[\mathcal J(y^*)]$ if $\mathcal R$ is Fréchet differentiable at $\mathcal J(y^*)$.

In addition, there exist
\[
p^*\in L^{\tilde p}(\Omega;V),\qquad
\mu^*\in L^{\tilde p}(\Omega;V^*),
\]
satisfying
\[ q^*  = \pi^*p^*, \qquad \lambda^* = \pi^*\mu^*,\]
such that
\begin{subequations}\label{eq:system_p_star}
    \begin{align}
    A^*p^*+\mu^*-J_y(y^*,\cdot) &=0&&\text{a.s.},\label{eq:adjHat}\\
\E[\pi^*\langle B^*p^*,u^*-v\rangle_{U^*,U}]
+\varrho'(u^*)(u^*-v)&\leq0 &&\forall v\in U_{ad},\label{eq:viHat}\\
\langle\mu^*,p^*\rangle&\geq0 &&\text{a.s.},\label{eq:signHat}\\
\langle\zeta^*,(p^*)^+\rangle
=\langle\zeta^*,(p^*)^-\rangle&=0 &&\text{a.s.},\label{eq:zetaHatPM}\\
\langle\mu^*,v\rangle&=0
&&\forall v\in V:\ v=0\ \text{q.e. on }\mathcal A^*, \text{ a.s.}\label{eq:exactSupportHat}
    \end{align}
\end{subequations}
\end{theorem}

\noindent Note that the condition \eqref{eq:lambdastar_orth_qw} should be understood as: for almost every $\omega \in \Omega$,
\[\langle\lambda^*(\omega),v\rangle=0 \qquad \forall v\in V:\ v=0\ \text{q.e. on }\mathcal A^*(\omega),\]%
and likewise for \eqref{eq:exactSupportHat}.

Let us make some initial observations.
\begin{itemize}[itemsep=0.1pt]
    \item In the system \eqref{eq:factorisedSystem}, we refer to $\pi^*$ as the \textit{risk identifier}, $q^*$ as the \textit{weighted adjoint}, and in \eqref{eq:system_p_star}, we call $p^*$  the \textit{unweighted adjoint} as it is a factor in $q^*$.

\item From  \eqref{eq:lambdastar_orth_qw} and \eqref{eq:exactSupportHat}, we obtain in particular
\[
\langle\lambda^*(\omega),y^*(\omega)-\psi(\omega)\rangle = \langle\mu^*(\omega),y^*(\omega)-\psi(\omega)\rangle
=0
\qquad\text{a.s.}
\]
\item The almost-sure C-sign conditions imply corresponding
inequalities in expectation whenever the duality
products are integrable. Indeed, for $q=2$,
H\"older's inequality yields
\[
\mathbb E[|\langle\lambda^*,q^*\rangle|]
\le
\|\lambda^*\|_{L^2(\Omega;V^*)}
\|q^*\|_{L^2(\Omega;V)}<\infty,
\]
and hence \eqref{eq:lambdaqgeq0Nonsmooth} holds also in expectation, i.e.,
\[
\mathbb E[\langle\lambda^*,q^*\rangle]\ge0.
\]
Similarly, if $\tilde p\ge2$, then
\[
\mathbb E[|\langle\mu^*,p^*\rangle|]
\le
\|\mu^*\|_{L^{\tilde p}(\Omega;V^*)}
\|p^*\|_{L^{\tilde p}(\Omega;V)}<\infty,
\]
so that the expected value version of \eqref{eq:signHat} is also a consequence:
\[
\mathbb E[\langle\mu^*,p^*\rangle]\ge0.
\]
The condition $\tilde p\ge2$ is equivalent to
$p\ge2$, or $1\le p<2$ and
$q\le 2p/(2-p)$. In particular, the restriction
$q\le2p$ or $p\ge2$ imposed in \cite{AGHS}, see \eqref{ass:furtherRestrictionOnPq} and the discussion below,
is sufficient but not necessary.

\end{itemize}
\begin{remark}
    A formal argument involving Clarke calculus applied to our problem reveals that \eqref{eq:factorisedSystem} is exactly the `right' system to expect as a stationarity system. Indeed, rewriting \eqref{eq:OCProblem} as
\[
\begin{aligned}
\min_{u,y,\zeta,X}\quad&
\cR[X]+\varrho(u)\\
\text{subject to}\quad&
Ay-f-Bu+\zeta=0,\\
&X-\cJ(y)=0,\\
&\zeta-(\zeta+c(y-\psi))^+=0,\\
&u\in U_{ad},
\end{aligned}
\]
for some $c>0$, and letting $\pi$, $q$ and $\eta$ denote the Lagrange
multipliers associated respectively with
the loss constraint, the state equation
and the complementarity equation, we can form the Lagrangian
\[
\begin{aligned}
L(u,y,\zeta,X;\pi,q,\eta)
={}&\cR[X]+\varrho(u) + \mathbf{1}_{U_{ad}}(u) +\E[\pi(\cJ(y)-X)]-\E[\langle q,Ay-f-Bu+\zeta\rangle]\\
&+\E\left[
\int_D\eta\,
\bigl(\zeta-(\zeta+c(y-\psi))^+\bigr)
\,\mathrm dx
\right],
\end{aligned}
\]
where $\mathbf{1}_{U_{ad}}$ is the indicator function of the admissible set. 
By considering stationarity (using standard calculations) we are led exactly to the system \eqref{eq:factorisedSystem}. Hence the system for $q^*$ is, in this sense, the natural C-stationarity system to expect.
\end{remark}

\subsubsection*{Comparison with literature and discussion of contributions}

As mentioned, this work builds upon what was done in \cite{AGHS}. We briefly compare our results to that paper. In \cite[Theorem 4.8]{AGHS}, the following so-called \textit{$\mathcal{E}$-almost weak stationarity} system was established for any local minimiser  {$u^* \in U_{ad}$ of \eqref{eq:OCProblem} with associated state $y^* \in  L^q(\Omega;V)$}: there exist  $\zeta^* \in L^q(\Omega;V^*)$, $q^* \in L^{q'}(\Omega;V)$,  $\pi^* \in L^{p'}(\Omega)$ and $\lambda^* \in L^{q'}(\Omega;V^*)$ such that 
\eqref{eq:stateCNew}---\eqref{eq:system_q_star_control_ineq} holds, but in lieu of our conditions \eqref{eq:lambdaqgeq0Nonsmooth}, \eqref{eq:zeta_q_pos_neg_zero} and \eqref{eq:lambdastar_orth_qw}, only the following were established:
\begin{subequations}\label{eq:eaWeakStationarity}
\setcounter{equation}{5}
\begin{align}
           \limsup_{\tau \to 0}\E[\langle \lambda_\tau^*, q_\tau^* \rangle] &\geq 0,\label{eq:limsupCond}\\
 \E [\langle\zeta^*, q^*\rangle]  &= 0\quad \text{if $q=2$,} \label{eq:zetaQ}\\
     \E[\langle \lambda^*, y^*-\psi \rangle] &= 0,\label{eq:new_label}\\
       \forall \epsilon > 0, \exists E^\epsilon \subset \mathcal{I}^* \text{ with } |\mathcal{I}^* \setminus E^\epsilon| \leq \epsilon : \E[\langle \lambda^*, v \rangle] &= 0 \quad \, \forall v \in L^{q}(\Omega;V)  :  v = 0 \text{ a.s.-a.e. on $\Omega\times D \setminus E^\epsilon$}.\label{eq:eaEaQ}
\end{align}
\end{subequations}
Above, $\mathcal{I}^*$ denotes the inactive/non-coincidence set associated to $y^*$. In addition, under the assumption that
\begin{equation}
   q \leq 2p \quad\text{or} \quad p \geq 2\label{ass:furtherRestrictionOnPq}
\end{equation}
\cite[Proposition 4.9]{AGHS} showed: there exists $\bar p \in L^{\tilde{p}}(\Omega;V)$  and $\bar\mu \in L^{\tilde p}(\Omega;V^*)$ such that
\begin{subequations}\label{eq:ecsSytem}
\begin{align}
    {A(\omega)^*}\bar{p}(\omega) + \bar{\mu}(\omega)  &= J_y(y^*(\omega),\omega) \quad \text{a.s.},\label{eq:adjoint-smooth}\\
        \E [\langle\zeta^*, \bar{p}\rangle] &= 0,\label{eq:zetaP}\\
       \E[ \langle \bar{\mu} , \bar{p} \rangle] &\geq 0,\label{eq:TT}\\
     \E[\langle \bar{\mu}, y^*-\psi \rangle] &= 0,\label{eq:complementarityLambda}\\
       \forall \epsilon > 0, \exists E^\epsilon \subset \mathcal{I}^* \text{ with } |\mathcal{I}^* \setminus E^\epsilon| \leq \epsilon : \E[\langle \bar{\mu}, v \rangle] &= 0 \, \forall v \in L^{\tilde p'}(\Omega;V)  : v = 0 \text{ a.s.-a.e. on $\Omega\times D \setminus E^\epsilon$}.\label{eq:eaEa}
\end{align}
\end{subequations}
which in combination with the already-stated system was classed as a \textit{weak $\mathcal{E}$-almost C-stationarity} system.  Furthermore,  under the assumptions that %
\begin{equation}
        \partial \cR[z] = \{\cR'[z]\} \quad  \forall z \in L^p(\Omega), \qquad \text{and}\qquad  \cR'\colon L^p(\Omega) \to L^{p'}(\Omega) \text{ is continuous} \label{ass:onRBeingNice}
\end{equation}
(and also %
\eqref{ass:furtherRestrictionOnPq}), it was shown  \cite[Proposition 4.11]{AGHS} that 
\[q^* = \pi^*\bar{p}, \qquad \lambda^* = \pi^*\bar{\mu},\]
providing a connection between \eqref{eq:eaWeakStationarity} and \eqref{eq:ecsSytem}.

Let us point out in what sense our work improves upon this.
\begin{enumerate}[label=(\alph*)]
    \item The sign condition
      $\E[ \langle \lambda^* , q^* \rangle] \geq 0$ 
was not shown, only \eqref{eq:limsupCond}. 

\emph{We obtain a stronger a.s. sign condition \eqref{eq:lambdaqgeq0Nonsmooth} on the limit.}
\item The orthogonality condition in \eqref{eq:zetaQ} holds only when $q=2$.  

\emph{We obtain an orthogonality condition which is much stronger: we prove an a.s. orthogonality for the positive and negative parts and we  do not need $q=2$ for the integrated condition ($q \in [2,\infty)$ suffices).}
\item The orthogonality conditions $\E[ \langle \zeta^*, (\bar{p})^+ \rangle] = \E[ \langle \zeta^*, (\bar{p})^- \rangle] = 0$ were not shown.

\emph{We obtain the missing conditions, in a stronger a.s. sense, see \eqref{eq:zetaHatPM}.}

\item The conditions \eqref{eq:eaEaQ} and \eqref{eq:eaEa} are $\mathcal{E}$-almost conditions that involve an unwanted $\epsilon > 0$ small set which arises as a result of using Egorov's theorem.

\emph{We remove the need for the $\mathcal{E}$-almost conditions and strengthen the statements substantially in \eqref{eq:lambdastar_orth_qw} and \eqref{eq:exactSupportHat}.}
\item  As pointed out in \cite{AGHS}, the assumption \eqref{ass:onRBeingNice}, used to deduce $q^*=\pi^*\bar{p}$, is quite strong as it essentially requires smoothness of $\cR$ everywhere. 

\emph{We weaken this requirement by needing only either $\cR$ to be Fréchet differentiable at $\cJ(y^*)$ (without needing \eqref{ass:onRBeingNice} or \eqref{ass:furtherRestrictionOnPq}), or needing $\cR$ to be a monotone risk measure.}
\end{enumerate}
\subsection{Examples}
Let us give an example setting to which our theory is applicable. Take $V=H^1_0(D)$ and let
$a\colon\Omega\times D\to\R$ satisfy
\[
a_0\leq a(\omega,x)\leq a_1 \qquad\text{for a.e.~$(\omega,x)\in\Omega\times D$,}
\]
where $0<a_0<a_1$ are constants.
Assume furthermore that
$\omega\mapsto a(\omega,\cdot)$ is strongly
measurable as an $L^\infty(D)$-valued map.
Define the operator $A(\omega)\colon V\to V^*$ by
\[
\langle A(\omega)y,z\rangle
=\int_D a(\omega,x)\grad y\cdot\grad z\D x,
\qquad \forall y,z\in V.
\]
Set $U=H=L^2(D)$ and consider the box-constrained
admissible set
\[
U_{ad}:=\{u\in L^2(D):
u_a\leq u\leq u_b\text{ a.e. in }D\},
\]
where $u_a,u_b\in L^2(D)$ satisfy
$u_a\leq u_b$ a.e. in $D$.
Take $B$ to be the canonical compact embedding
$L^2(D)\ctsCompact H^{-1}(D)$, i.e.,
$B(\omega)u\equiv u$ viewed as an element of
$V^*=H^{-1}(D)$.

Let
\[
f\in L^2(\Omega;H^{-1}(D)),\qquad
\psi\in L^2(\Omega;H^1_0(D)),
\]
and choose the exponents $q=2$ and $p=1$.
Define the tracking-type functional and regulariser
\begin{equation*}
J(y):=\frac12\norm{y-y_d}{H}^2,
\qquad
\varrho(u):=\frac{\nu}{2}\norm{u}{H}^2,
\end{equation*}
where $y_d\in L^2(D)$ is a given desired state
and $\nu>0$ is a constant. The assumptions on the operators, objective
and control regularisation are satisfied.
Moreover, $U_{ad}$ is nonempty, closed,
convex and bounded.

We now give several examples of risk measures $\cR\colon L^1(\Omega)\to\R$ to which our results apply.
\begin{itemize}[leftmargin=*]
\item \textbf{Expectation.}
The risk-neutral choice 
\[\cR[X]=\E[X]\]
is linear, monotone and continuously Fr\'echet differentiable.
Thus \cref{thmmain} applies with both the differentiability and the monotonicity assumption satisfied.

\item \textbf{Conditional value-at-risk (CVaR).}
For $\beta\in[0,1)$, define
\begin{equation*}%
\cR[X]=\CVaR_\beta[X]
:=\inf_{t\in\R}\left\{t+\frac{1}{1-\beta}\E[(X-t)^+]\right\}.
\end{equation*}
CVaR is finite, convex, monotone and Lipschitz continuous on $L^1(\Omega)$; see
\cite[\S 6.2.4]{Shapiro2009}. For $\beta=0$ it reduces to the expectation.
Thus \cref{thmmain} applies under monotonicity of $\cR$.

\item \textbf{Finite mixtures of CVaRs and mean--CVaR.}
For $\beta_i\in[0,1)$ and weights $w_i\ge0$ with $\sum_{i=1}^m w_i=1$, let
\[
\cR[X]=\sum_{i=1}^m w_i\CVaR_{\beta_i}[X].
\]
This is again finite, convex, monotone and continuous on $L^1(\Omega)$, so \cref{thmmain}
applies under monotonicity of $\cR$. In particular, this includes the mean--CVaR combinations
\begin{equation}\label{eq:appMeanCVaR}
\cR_{\alpha,\beta}(X):=\alpha\E[X]+(1-\alpha)\CVaR_\beta[X],
\qquad\alpha\in[0,1],\ \beta\in(0,1),
\end{equation}
see \cite[Example 6.16]{Shapiro2009} and \cref{sec:mean_CVaR}.

\item \textbf{Spectral risk measures.}
Let $\phi\colon[0,1)\to\R_+$ be bounded and non-decreasing with $\int_0^1\phi(p)\,\D p=1$
(i.e., $\phi$ is a \emph{risk spectrum}), and let $\mathrm{VaR}_p[X]:=\inf\{t\in\R:\mathbb P(X\le t)\ge p\}$. Then
\[
\cR[X]=\int_0^1\phi(p)\,\mathrm{VaR}_p[X]\;\mathrm{d}p
\]
is a finite, convex, monotone and continuous risk measure on $L^1(\Omega)$; see \cite{Acerbi}.
Hence \cref{thmmain} applies under monotonicity of $\cR$.

\item \textbf{Mean--upper-semideviation of order one.}
For $c\in[0,1]$, let
\[
\cR[X]=\E[X]+c\,\E\big[(X-\E[X])^+\big].
\]
This is finite, convex, monotone and Lipschitz continuous on $L^1(\Omega)$ (see
\cite[\S 6.2.2]{Shapiro2009}), so \cref{thmmain} applies under monotonicity of $\cR$.

\item \textbf{Finite maxima of expectations.}
Let $Z_1,\ldots,Z_m\in L^\infty(\Omega)$ be nonnegative with $\E[Z_i]=1$, and let
\[
\cR[X]=\max_{1\le i\le m}\E[Z_iX].
\]
This is finite, convex, monotone and Lipschitz continuous on $L^1(\Omega)$, so \cref{thmmain}
applies under monotonicity of $\cR$. Moreover, if the maximum of $i\mapsto\E[Z_i\cJ(y^*)]$ is
attained at a unique index $i^*$, then $\cR$ is Fr\'echet differentiable at $\cJ(y^*)$. In that case \cref{thmmain} also applies under
differentiability of $\cR$.

\item \textbf{Log-sum-exp risk measure.}
Let $Z_1,\ldots,Z_m\in L^\infty(\Omega)$ be nonnegative with $\E[Z_i]=1$, let $w_i>0$ with
$\sum_{i=1}^m w_i=1$, and let $\gamma>0$. Define
\[
\cR_\gamma(X)=\frac1\gamma\log\Big(\sum_{i=1}^m w_i\exp\big(\gamma\,\E[Z_iX]\big)\Big).
\]
This risk measure is finite, convex, monotone and continuously Fr\'echet differentiable on
$L^1(\Omega)$, so both the differentiability and the monotonicity assumption of \cref{thmmain}
are satisfied. On a finite probability space $\Omega=\{\omega_1,\ldots,\omega_m\}$ with
$\mathbb P(\omega_i)>0$, the choice $Z_i=\mathbf 1_{\{\omega_i\}}/\mathbb P(\omega_i)$ and
$w_i=\mathbb P(\omega_i)$ gives $\E[Z_iX]=X(\omega_i)$ and recovers the \textbf{entropic risk measure}
$\frac1\gamma\log\E[\exp(\gamma X)]$.

\item \textbf{General class of worst-case expectations.}
Some of the examples above are instances of \emph{worst-case expectations}. Let
$\mathcal Q\subset L^\infty(\Omega)$ be nonempty and uniformly bounded with
\[
Z\ge0,\qquad\E[Z]=1\qquad\text{for all }Z\in\mathcal Q,
\]
and define
\[
\cR[X]=\sup_{Z\in\mathcal Q}\E[ZX].
\]
Then $\cR$ is finite, convex, monotone and Lipschitz continuous on $L^1(\Omega)$, %
so \cref{thmmain} applies under monotonicity of $\cR$.

\item \textbf{Conditional-mean penalty.}
Let $B\in\mathcal F$ with $0<\mathbb P(B)<1$, set $Z_B:=\mathbf 1_B/\mathbb P(B)$ and let $c>0$. Define
\[
\cR[X]=\E[X]+c\big(\E[Z_BX]-\E[X]\big)^2 .
\]
This is non-monotone, finite, convex and continuously Fr\'echet differentiable on $L^1(\Omega)$ with
\[
\cR'[X]=1+2c\,\big(\E[Z_BX]-\E[X]\big)(Z_B-1).
\]
In particular, the risk identifier
$\pi^*=\cR'[\cJ(y^*)]$ may take negative values, so \cref{thmmain} applies only under
differentiability of $\cR$.

\end{itemize}
Thus our C-stationarity result
applies to a broad class of convex  risk measures.%

\subsection{Mean--CVaR control of a random obstacle problem}\label{sec:mean_CVaR}
Mean--CVaR risk measures of the form \eqref{eq:appMeanCVaR} are among the most widely used in practice: through the parameter $\alpha$ they interpolate between risk-neutral optimisation and control of the tail of the loss distribution, while retaining a clear interpretation in terms of the average of the worst $(1-\beta)$-fraction of outcomes. Mean--CVaR appears in application areas such as finance \cite{krokhmal2002portfolio}, energy systems \cite[\S 4.2.2]{MR3236311}, and machine learning \cite[\S C]{Chow2015}.

We examine \eqref{eq:OCProblem} with mean--CVaR chosen as the risk measure. We apply \cref{thmmain} to the elliptic obstacle problem and the quadratic tracking objective introduced above. In particular, $V=H^1_0(D)$, $U=L^2(D)$, $B$ is the canonical embedding into $V^*$, and
\[
\langle A(\omega)y,v\rangle
 =\int_D a(\omega,x)\nabla y\cdot\nabla v\,\mathrm dx,
\qquad
J(y)=\frac12\|y-y_d\|_{L^2(D)}^2,
\]
\[
\varrho(u)=\frac\nu2\|u\|_{L^2(D)}^2, \qquad 
U_{ad}=\{u\in L^2(D):u_a\le u\le u_b\text{ a.e.}\}.
\]
Let $\beta\in(0,1)$ and $\alpha\in[0,1]$, and recall the mean--CVaR risk measure \eqref{eq:appMeanCVaR}:
\begin{equation*}
\cR_{\alpha,\beta}(X)
 :=\alpha\E[X]+(1-\alpha)\CVaR_\beta[X].
\end{equation*}
Thus $\alpha=1$ gives the risk-neutral expectation, whereas
$\alpha=0$ gives pure CVaR. As it is finite, convex,
and monotone on $L^1(\Omega)$, \cref{thmmain} applies
with $p=1$, $q=2$ and $\tilde p=2$. Let $u^*\in U_{ad}$ be a local minimiser of \eqref{eq:OCProblem} with associated state $y^*=S(u^*)$. 

Recall that, for $\beta\in(0,1)$, a
\textit{$\beta$-quantile} of a random variable $X$ is any threshold $t$ where at most a $\beta$-fraction of the probability lies strictly below $t$, and at least a $\beta$-fraction lies at or below $t$. In other words, it is any $t\in\R$ satisfying
\[
\pP(X<t)\leq\beta\leq\pP(X\leq t)\qquad\text{or equivalently}\qquad
\pP(X>t)\leq1-\beta\leq\pP(X\geq t).
\]
In $\CVaR_\beta[\cJ(y^*)]$, the parameter $\beta\in(0,1)$ specifies the
proportion $1-\beta$ of the largest tracking
losses considered by CVaR. The corresponding
loss threshold is given by a $\beta$-quantile
$t^*\in\R$ of $\cJ(y^*)$, characterised by
\[
\pP(\cJ(y^*)<t^*)\leq\beta
\leq\pP(\cJ(y^*)\leq t^*).
\]
In particular, scenarios with
$\cJ(y^*)>t^*$ belong to the upper loss tail. CVaR represents the average loss over the
worst $1-\beta$ fraction of the probability
distribution.

If $\pP(\cJ(y^*)=t^*)=0$, then
\[
\CVaR_\beta[\cJ(y^*)]
=
\E[\cJ(y^*)\mid \cJ(y^*)>t^*].
\]
If instead $\pP(\cJ(y^*)=t^*)>0$, the loss
distribution is said to have an \emph{atom}
at $t^*$, meaning that the tracking loss
equals the threshold with positive probability.

Applying \cref{thmmain} in this setting, we can prove the following properties on the optimal control and risk identifier.
\begin{prop}
\label{prop:appMeanCVaR}
The multipliers in \cref{thmmain}
may be chosen together with a
$\beta$-quantile $t^*\in\R$ of $\cJ(y^*)$
and $\xi^*\in L^\infty(\Omega)$ such that
the following properties hold. %

\medskip
\noindent The risk multiplier satisfies
\begin{equation}\label{eq:appRiskDensity}
\pi^*=\alpha+(1-\alpha)\xi^*,
\end{equation}
where
\begin{equation*}%
\xi^*(\omega) \in
\begin{cases}
\{0\}, &\text{if } \cJ(y^*(\omega))<t^*,\\
[0,(1-\beta)^{-1}],
  &\text{if } \cJ(y^*(\omega))=t^*,\\
  \{(1-\beta)^{-1}\},&\text{if } \cJ(y^*(\omega))>t^*,
\end{cases}
\qquad\text{and}\qquad \E[\xi^*]=1.
\end{equation*}
Thus, the following complementarity conditions hold a.s.:
\begin{equation}\label{eq:appRiskComplementarity}
\begin{aligned}
\alpha\leq\pi^*&\leq \alpha+\frac{1-\alpha}{1-\beta},\\
(\pi^*-\alpha)(t^*-\cJ(y^*))^+&=0,\\
\left(\alpha+\frac{1-\alpha}{1-\beta}-\pi^*\right)(\cJ(y^*)-t^*)^+&=0,\\
\E[\pi^*]&=1.
\end{aligned}
\end{equation}
In addition, the optimal control has the pointwise representation
\begin{equation}\label{eq:appProjection}
u^*(x)=\operatorname{proj}_{[u_a(x),u_b(x)]}
\left(-\frac1\nu\E[\pi^*p^*](x)\right)
\quad\text{for a.e. }x\in D.
\end{equation}
The quantity entering the projection is the mixture
\[
\E[\pi^*p^*]
=\alpha\E[p^*]+(1-\alpha)\E[\xi^*p^*],
\]
which, if $\pP(\cJ(y^*)=t^*)=0$, reads
\begin{equation}\label{eq:appMeanTail}
\E[\pi^*p^*]
 =\alpha\E[p^*]+(1-\alpha)\E[p^*\mid \cJ(y^*)>t^*],
\end{equation}
so that in this case, \eqref{eq:appProjection} averages the almost-sure adjoints
using a mixture of their ordinary expectation and their
conditional expectation over high-loss scenarios.
\end{prop}

\begin{proof}
We begin with
\[
\partial\cR_{\alpha,\beta}[X]
=\alpha\mathbbm1
+(1-\alpha)\partial\CVaR_\beta(X).
\]
Recall the dual representation
\[
\CVaR_\beta[X]
=\sup\left\{\E[\xi X]:
\xi\in L^\infty(\Omega), \;
0\leq\xi\leq(1-\beta)^{-1},\; \E[\xi]=1\right\}
\]
and that the
subdifferential of $\CVaR_\beta$ can be written as
\[
\begin{aligned}
\partial\CVaR_\beta[X]
=\big\{\xi\in L^\infty(\Omega):\;&
0\leq\xi\leq(1-\beta)^{-1},\
\E[\xi]=1,\\
&\xi=0\quad\text{on }\{X<t\},\\
&\xi=(1-\beta)^{-1}
\quad\text{on }\{X>t\}\big\}
\end{aligned}
\]
for any $\beta$-quantile $t$ of $X$, see \cite[Example 6.16]{Shapiro2009}.
Applying this characterisation with $X=\cJ(y^*)$
gives \eqref{eq:appRiskDensity}--%
\eqref{eq:appRiskComplementarity}.

Since $B$ is the
canonical embedding,
the control variational inequality becomes
\[
(\nu u^*+\E[\pi^*p^*],v-u^*)_{L^2(D)}\geq0
\qquad\forall v\in U_{ad},
\]
which is equivalent to the box projection
\eqref{eq:appProjection}.
If $\pP(\cJ(y^*)=t^*)=0$, the quantile relation implies
$\pP(\cJ(y^*)>t^*)=1-\beta$, and therefore
$\xi^*=(1-\beta)^{-1}\mathbbm1_{\{\cJ(y^*)>t^*\}}$
almost surely. Formula \eqref{eq:appMeanTail} then follows.
\end{proof}
Let us make a few observations:
\begin{itemize}[itemsep=0.1pt]
    \item 
    Since $\pi^*\geq0$ a.s. and
$\E[\pi^*]=1$, it defines a probability
measure $\mathbb Q^*$ on $(\Omega,\mathcal F)$ by
\[
\mathbb Q^*(A)
:=\E[\pi^*\mathbf1_A],
\qquad A\in\mathcal F.
\]
In particular, $\mathbb Q^*\ll\pP$ and
\[
\frac{\mathrm d\mathbb Q^*}{\mathrm d\pP}=\pi^*.
\]
The pointwise control formula \eqref{eq:appProjection} can be written as
\[
u^*=\operatorname{proj}_{[u_a,u_b]}
\left(-\frac1\nu\E_{\mathbb Q^*}[p^*]\right).
\]
When $\pP(\cJ(y^*)=t^*)=0$,
\[
\mathbb Q^*=\alpha\pP
 +(1-\alpha)\pP(\,\cdot\mid \cJ(y^*)>t^*).
\]
In particular, the risk-neutral case $\alpha=1$ uses
$\E[p^*]$, whereas the pure CVaR  case $\alpha=0$ uses
$\E[p^*\mid \cJ(y^*)>t^*]$. 

\item The inequalities \eqref{eq:appRiskComplementarity} describe complementarity in probability space: scenarios with losses below the threshold defining
the worst $1-\beta$ fraction receive the baseline
weight $\alpha$, whereas scenarios with losses
strictly above this threshold  receive the larger weight $\alpha+\frac{1-\alpha}{1-\beta}$.  
Scenarios with losses exactly at the threshold
may receive intermediate weights.

\item

The risk identifier $\pi^*$ determines how the
  adjoints are weighted in the control
optimality condition. 

When $\alpha>0$, we have $\pi^*\geq\alpha>0$, so the unweighted adjoint 
can be recovered from the weighted adjoint through
$p^*=q^*/\pi^*$. In contrast, for pure CVaR
($\alpha=0$), scenarios with tracking losses
strictly below the threshold $t^*$ receive zero
weight, i.e., $\pi^*=0$ on $\{\cJ(y^*)<t^*\}$,
and thus $q^*=\lambda^*=0$ there.
Although these scenarios do not contribute
directly to the control optimality condition,
\cref{thmmain} still provides an unweighted
adjoint $p^*$ and multiplier $\mu^*$ satisfying
the a.s. adjoint equation and C-stationarity
conditions. 

This illustrates the importance
of constructing the unweighted adjoint on
$\{\pi^*=0\}$, where it cannot be recovered
by division from the weighted adjoint.
\end{itemize}

\begin{remark}[Scenarios at the CVaR threshold]
\label{rem:appQuantileAtoms}
When $\pP(\cJ(y^*)=t^*)>0$, the risk multiplier on the quantile
level set may be nonunique. More precisely,
\begin{equation*}%
\begin{aligned}
\E[\pi^*p^*]
={}&\alpha\E[p^*]
 +\frac{1-\alpha}{1-\beta}
   \E[p^*\mathbbm1_{\{\cJ(y^*)>t^*\}}] +\E[(\pi^*-\alpha)p^*
  \mathbbm1_{\{\cJ(y^*)=t^*\}}],
\end{aligned}
\end{equation*}
and normalisation fixes the total additional weight:
\begin{equation*}%
\E[(\pi^*-\alpha)\mathbbm1_{\{\cJ(y^*)=t^*\}}]
=(1-\alpha)\left(
1-\frac{\pP(\cJ(y^*)>t^*)}{1-\beta}\right).
\end{equation*}
The values of $\pi^*$ on the quantile level set $\{\mathcal J(y^*)=t^*\}$ need not be uniquely determined, but must satisfy the coupled C-stationarity conditions.

For example, consider four equally likely scenarios
and pure CVaR with $\beta=1/2$, so that CVaR
averages the worst $50\%$ of the losses. Suppose
\[
\cJ(y_1^*)<\cJ(y_2^*)=\cJ(y_3^*)<\cJ(y_4^*).
\]
The largest loss occurs in scenario $4$, while
scenarios $2$ and $3$ have identical losses at
the threshold $t^*$. Scenario $4$ accounts for $25\%$ of the probability mass, so another $25\%$ is needed to reach the worst $50\%$ represented by CVaR. Since scenarios $2$ and $3$ have identical losses at the threshold, this additional probability mass can be distributed between them without changing the CVaR value. 
This freedom is reflected in the choice of the risk multiplier $\pi^* \in\partial\mathcal R[\mathcal J(y^*)]$, whose components specify the scenario weights. More precisely, the admissible risk subgradients are given by
\[
(\pi_1^*, \pi_2^*, \pi_3^*,\pi_4^*)
=(0,\theta,2-\theta,2),
\qquad 0\leq\theta\leq2.
\]
Here $\theta$ determines the allocation of risk weight between scenarios $2$ and $3$.
In particular, $\theta=1$
assigns equal weights, whereas $\theta=0$
or $\theta=2$ assigns all the remaining weight
to one of the scenarios.

The CVaR subdifferential permits any such
allocation, but the control optimality condition
requires a choice compatible with the associated
adjoints. Indeed,
\[
\E[\pi^*p^*]
=\frac14\bigl(
\theta p_2^*+(2-\theta)p_3^*+2p_4^*
\bigr).
\]
This shows that the CVaR
subdifferential may contain multiple risk
multipliers when scenarios have identical losses
at the quantile threshold. The particular
multiplier must be compatible with the coupled
C-stationarity system. 
\end{remark}

\section{Preliminaries: existence and first stationarity results}\label{sec:optimisation-problem}
Let us proceed on the road to proving our theorem. We will in this section define and study basic properties of an approximated problem to \eqref{eq:OCProblem}. First, we discuss the well posedness and properties of \eqref{eq:OCProblem} itself.

We denote the solution mapping $u \mapsto y(\omega)$ in \eqref{eq:VIProblemFormal} by $S_\omega \colon U \to V$ and its associated superposition operator $S$ by
\begin{align*}
    S(u)(\omega) := S_\omega(u).    %
\end{align*}
Note that this is well defined by the standard Lions--Stampacchia theory. In \cite[\S 2.1]{AGHS}, the following basic properties of $S$ are shown.
\begin{prop}\label{prop:propertiesOfS}
The solution map $S$ %
is well defined and satisfies the following properties.
\begin{enumerate}[label=(\roman*)]\itemsep0em 
    \item  $S(u) \in L^{\min(r,s)}(\Omega;V)$ for all $u \in U$.
    \item The following estimate holds on an a.s. level:
    \[\norm{S_\omega(u) - S_\omega(v)}{V} 
    \leq C_a^{-1} \norm{B(\omega) (u-v)}{V^*}.\] %
    \item   
    If $u_n \weaklyto u$ in $U$, then 
    \[S_\omega(u_n) \to S_\omega(u)\text{ in $V$ a.s.}\]
    and
    \begin{align*}
    S(u_n)  &\to S(u) \text{ in $L^{q}(\Omega;V)$},\qquad   S(u_n)-S(u)  \to 0\text{ in $L^t(\Omega;V)$ {for all $t \in [1, \infty)$}}.
    \end{align*}
\end{enumerate}
\end{prop} 
The existence of optimal controls is shown in \cite[\S 2]{AGHS} with the direct method. We record the result here. 
\begin{prop}\label{prop:existenceOC} %
Assume that either $U_{ad}$ is bounded or $\mathcal{R} \circ \cJ \circ S  + \varrho $ is coercive. Then %
there exists an optimal control $u^*$ of \eqref{eq:OCProblem}.
\end{prop}
The coercivity part of the assumption in the above existence result %
is automatically met in the setting where $\cR$ is monotonic with $\cR[0] \geq 0$, $J$ is a tracking-type functional %
and $\varrho(u) = (\nu\slash 2)\norm{u}{U}^2$ for some constant $\nu>0$.
 
\subsection{Penalised problem and stationarity} \label{sec:penalisedProblemStat}
As in \cite[\S 4]{AGHS}, our approach to deriving stationarity conditions for minimisers of the problem \eqref{eq:OCProblem} is through penalising the variational inequality via semilinear PDEs, establishing the penalised KKT system, and then passing to the limit. Note that this general idea of a penalty procedure in the deterministic case is classical, see e.g. \cite[\S 3.5.2, p. 370]{Lions1969}, \cite[Chapter 1, \S 3.2]{GlowinskiLT} (as well as the recent abstract framework \cite{AW}). 

In order to penalise the inequality constraint in \eqref{eq:VIProblemFormal}, we define the following differentiable approximation of the positive part  function $(\cdot)^+ \equiv \max(0,\cdot)$:
\begin{equation*}%
m_\tau(r) %
:= \begin{cases}
0 &: r \leq 0,\\
\frac{r^2}{2\tau} &: 0 < r < \tau,\\
r-\frac{\tau}{2} &: r \geq \tau.
\end{cases}
\end{equation*}
We have $m_\tau \in C^1(\mathbb{R})$, $m_\tau' \in [0,1]$ and   $m_\tau\colon H^1(D) \to L^2(D)$ is $C^1$.  Other choices of penalisations or smooth approximations $m_\tau$ are possible, see \cite{AW}.

Then we consider for a fixed $\omega \in \Omega$ the penalised problem 
 \begin{equation}
 \label{eq:penalisedPDE}
 A(\omega)y +  \frac 1\tau m_\tau(y-\psi(\omega)) =f(\omega)+ B(\omega)u
 \end{equation}
 and denote the solution map $u \mapsto y$ as $T_{\tau,\omega} \colon U \to V$ and the associated superposition map %
 $T_\tau(u)(\omega) := T_{\tau,\omega}(u).$  The solution map is well defined since \eqref{eq:penalisedPDE} has a unique solution by \cite[Theorem 2.6]{Roubicek}.

Analogously to \eqref{eq:OCProblem}, we consider the regularised problem
 \begin{equation}\label{eq:penalisedOCProblem}
 \min_{u \in U_{ad}} \cR[\cJ(T_\tau(u))] +\varrho(u).
 \end{equation}
To guarantee the existence and boundedness of minimising sequences, let us assume that either $U_{ad}$ is bounded or that the family
$\{\mathcal R\circ\mathcal J\circ T_\tau+\varrho\}_{\tau>0}$
is coercive, with coercivity of the latter
uniform in $\tau$ for $\tau$ sufficiently small. We denote the optimal controls of the regularised problem \eqref{eq:penalisedOCProblem} by $u_\tau^* \in U_{ad}$,  with associated states 
 \[y_\tau^* := T_\tau(u_\tau^*) \in L^{\min(r,s)}(\Omega;V)\]
satisfying
 \begin{equation*}
      A(\omega)y_\tau^*(\omega) +  \frac 1\tau m_\tau(y_\tau^*(\omega)-\psi(\omega)) =f(\omega)+ B(\omega)u_\tau^*\quad \text{a.s.}%
 \end{equation*}
 See \cite[\S 3]{AGHS} for more details.  Combining Proposition 3.12 and Proposition 3.13 of \cite{AGHS}, we obtain the following KKT system for \eqref{eq:penalisedOCProblem}.

\begin{prop}[Penalised stationarity system]\label{prop:KKT-penalised}
There exist
\[
q_\tau^*\in L^{q'}(\Omega;V),\qquad \bar{p}_\tau\in L^{\tilde p}(\Omega;V),\qquad
\pi_\tau^*\in L^{p'}(\Omega),
\]
satisfying
\[q_\tau^* = \bar p_\tau\pi_\tau^*\]
such that
\begin{subequations}\label{eq:tauKKT}
\begin{align}
A(\omega)^*q_\tau^*(\omega) + \frac 1\tau m_\tau'(y_\tau^*(\omega)-\psi(\omega))q_\tau^*(\omega) &= J_y(y_\tau^*(\omega),\omega)\pi_\tau^*(\omega) && \text{a.s.},\label{eq:ptau}\\
 \pi_\tau^* &\in \partial \cR[\cJ(y_\tau^*)],\label{eq:tauSS4}\\
A(\omega)^*\bar{p}_\tau(\omega)+\frac1\tau m_\tau'(y_\tau^*(\omega)-\psi(\omega))\bar{p}_\tau(\omega)
&=J_y(y_\tau^*(\omega),\omega) &&\text{a.s.},\label{eq:ptau-old}\\
\E[\pi_\tau^*\langle B^*\bar{p}_\tau,u_\tau^*-v\rangle_{U^*,U}]
+\varrho'(u_\tau^*)(u_\tau^*-v)&\leq0\qquad\forall v\in U_{ad}.\label{eq:tauSS3}
\end{align}
\end{subequations}
\end{prop}
Regarding the limiting convergence behaviour as $\tau \downarrow 0$, we have for the primal variables by using weak compactness and a standard consistency argument the following \cite[Lemma 4.1]{AGHS}:
\begin{lemma}\label{lem:boundednessofyandutaus}
There exists $(y^*, u^*) \in L^q(\Omega; V) \times U$ such that (for a subsequence that we do not distinguish)
\begin{equation*}
    \begin{aligned}
    y_\tau^* &\to y^* &&\text{in $L^{q}(\Omega;V)$},\\
    u^*_\tau &\weaklyto u^* &&\text{in $U$},
    \end{aligned}
\end{equation*}
as $\tau \downarrow 0$, where $u^* \in U$ is a minimiser of \eqref{eq:OCProblem} and $y^* \in L^q(\Omega;V)$ is the associated state, i.e., $y^* = S(u^*)$.

\end{lemma}
The limiting behaviour for dual variables is as follows \cite[Lemma 4.4]{AGHS}.
\begin{lemma}
\label{lem:bounded-dual-variables}
There exist
\[
\bar{p}\in L^{\tilde p}(\Omega;V),\qquad
q^*\in L^{q'}(\Omega;V),\qquad
\pi^*\in L^{p'}(\Omega),
\]
such that (for a subsequence that we do not distinguish)
\[
\begin{aligned}
\bar{p}_\tau&\weaklyto \bar{p} &&\text{in }L^{\tilde p}(\Omega;V),\\
q_\tau^*&\weaklyto q^* &&\text{in }L^{q'}(\Omega;V),\\
\pi_\tau^*&\weaklyto\pi^* &&\text{in }L^{p'}(\Omega)
\quad(\weakstar\text{ if }p'=\infty)
\end{aligned}
\]
as $\tau\downarrow0$. 
\end{lemma}
Note that \eqref{eq:tauSS4} means 
\[\cR[g]-\cR[\cJ(y_\tau^*)]-\E[\pi_\tau^*(g-\cJ(y_\tau^*))]\geq0
\qquad\forall g\in L^p(\Omega),\]
and, using the above convergences, it is easy to pass to the limit here to obtain $\pi^* \in \partial \cR[\cJ(y^*)]$.

Now, let us define the adjoint multipliers
\begin{align}
\zeta_\tau^*(\omega) &:=\frac1\tau m_\tau(y_\tau^*(\omega)-\psi(\omega))=f(\omega)+Bu_\tau^*-Ay_\tau^*(\omega),\\
\lambda_\tau^*(\omega) &:= \frac 1\tau m_\tau'(y_\tau^*(\omega)-\psi(\omega))q_\tau^*(\omega) = J_y(y_\tau^*(\omega),\omega)\pi_\tau^*(\omega) - A(\omega)^*q_\tau^*(\omega),\\
\bar{\mu}_\tau(\omega)
&:=\frac1\tau m_\tau'(y_\tau^*(\omega)-\psi(\omega))\bar{p}_\tau(\omega)
=J_y(y_\tau^*(\omega),\omega)-A(\omega)^*\bar{p}_\tau(\omega),    \label{eq:lambdaTauDef}
\end{align}
and note
\[\lambda_\tau^* = \pi_\tau^*\bar{\mu}_\tau.\]
By the above convergences and the fact that $B$ is completely continuous, we immediately obtain
\begin{align}
\zeta_\tau^*\to\zeta^* &:= f+Bu^*-Ay^*
\qquad\text{ in }L^q(\Omega;V^*)\label{eq:zetaStrongNew}\\
\lambda_\tau^*\weaklyto
\lambda^* &:= \pi^*J_y(y^*,\cdot)-A^*q^*
\quad\text{in }L^{q'}(\Omega;V^*),\label{eq:weightedLambdaWeakLimit}\\
\bar{\mu}_\tau\weaklyto \bar{\mu} &:= J_y(y^*,\cdot)-A^*\bar{p} 
\qquad\text{ \;in }L^{\tilde p}(\Omega;V^*),\label{eq:lambdaWeakLimits}
\end{align}
along the selected subsequence.
In \cite[Assumption 3.8]{AGHS}, it was additionally assumed that, for every
$v\in U_{ad}$, the map
$u\mapsto\varrho'(u)(u-v)$ is weakly lower semicontinuous. This assumption was used to obtain \eqref{eq:system_q_star_control_ineq} by  passing to the limit in the VI relating the control to $q_\tau^*$ \eqref{eq:tauSS3}:  %
\[
\E[\langle B^*q_\tau^*,u_\tau^*-v\rangle_{U^*,U}]
+\varrho'(u_\tau^*)(u_\tau^*-v)\le0
\qquad\forall v\in U_{ad},
\]
see the proof of Theorem 4.8 of \cite{AGHS}.  
In fact, it is not needed: the already-assumed convexity and G\^ateaux differentiability 
on $\varrho$ suffice, as the next result shows.
\begin{prop}\label{prop:new_1}
We have
\[
\E[\langle B^*q^*,u^*-v\rangle_{U^*,U}]
+\varrho'(u^*)(u^*-v)\le0
\qquad\forall v\in U_{ad}.\]
\end{prop}
\begin{proof}
By convexity of $\varrho$, we have
\[
\varrho(u_\tau^*)+\varrho'(u_\tau^*)(v-u_\tau^*) \leq \varrho(v),
\]
and hence by \eqref{eq:tauSS3}
\begin{equation}\label{eq:eq1}
\varrho(u_\tau^*)
+\E[\langle B^*q_\tau^*,u_\tau^*-v\rangle_{U^*,U}] \leq \varrho(v) \quad\forall v\in U_{ad}.
\end{equation}
Since $u_\tau^*\rightharpoonup u^*$ in $U$, the complete continuity assumption
\eqref{ass:BcompletelyCts} combined with the uniform bounds in \cref{ass:newCombinedOnAandB} and a Dominated Convergence Theorem argument implies
\[
B(u_\tau^*-v)\to B(u^*-v)
\quad\text{in }L^q(\Omega;V^*).
\]
Together with $q_\tau^*\rightharpoonup q^*$ in
$L^{q'}(\Omega;V)$, this yields
\[
\E[\langle B^*q_\tau^*,u_\tau^*-v\rangle_{U^*,U}]
\to
\E[\langle B^*q^*,u^*-v\rangle_{U^*,U}].
\]
Furthermore, the convexity and continuity of $\varrho$ implies its weak
lower semicontinuity. Hence taking the limit inferior in the inequality \eqref{eq:eq1}, we obtain
\[
\varrho(u^*)
+\E[\langle B^*q^*,u^*-v\rangle_{U^*,U}] \leq \varrho(v)
\quad\forall v\in U_{ad}.
\]
Now fix $w\in U_{ad}$ and substitute
$v=u^*+t(w-u^*)\in U_{ad}$, where $t\in(0,1]$.
Dividing the resulting inequality by $t$ and letting
$t\downarrow0$, the G\^ateaux differentiability of
$\varrho$ gives
\[
\E[\langle B^*q^*,u^*-w\rangle_{U^*,U}]
+\varrho'(u^*)(u^*-w)\le0
\quad\forall w\in U_{ad},
\]
which is exactly the desired claim (which is \eqref{eq:system_q_star_control_ineq}).
\end{proof}
So far, we have shown conditions \eqref{eq:stateCNew}--\eqref{eq:system_q_star_control_ineq} (bearing in mind that the first two conditions are simply the complementarity form of the VI \eqref{eq:VIProblemFormal}). In the next sections we will prove the remaining statements.

 \section{Strengthened stationarity conditions}\label{sec:strengthened}

Now we depart from \cite{AGHS} in a more substantial manner.  Throughout this section, all limits as $\tau\downarrow0$
are understood along a fixed subsequence (which we do not relabel),
for which the convergences established in
\cref{sec:optimisation-problem} hold.
After further extraction, we may assume that
\[
y_\tau^*(\omega)\to y^*(\omega)
\quad\text{in }V,
\qquad
J_y(y_\tau^*(\omega),\omega)
\to J_y(y^*(\omega),\omega)
\quad\text{in }V^*,
\]
and
\[
\zeta_\tau^*(\omega)\to\zeta^*(\omega)
\quad\text{in }V^*
\]
hold almost surely. All subsequent pointwise arguments are
understood outside a common null set and further subsequences
will not be relabelled.

First, we will improve the conditions  \eqref{eq:zetaQ} and \eqref{eq:zetaP} obtained in \cite{AGHS} to an a.s. condition involving the absolute value of $\bar{p}$.  This resolves the issue raised in \cite[Remark 4.10]{AGHS}, and we even prove a stronger statement than discussed there as we show a.s. conditions, not only in expectation.
\begin{prop} \label{prop:separateContactOrthogonality}
We have
\begin{align}
\langle\zeta^*,|\bar{p}|\rangle&=0
\qquad\text{a.s.},\label{eq:pAbsOrthogonality}\\
\langle\zeta^*,|q^*|\rangle&=0
\qquad\text{a.s.},\label{eq:qAbsOrthogonality}
\end{align}
and hence
\begin{align*}
\langle\zeta^*,(\bar{p})^+\rangle
&=\langle\zeta^*,(\bar{p})^-\rangle=0
\qquad\text{a.s.},\\
\langle\zeta^*,(q^*)^+\rangle
&=\langle\zeta^*,(q^*)^-\rangle=0
\qquad\text{a.s.}    
\end{align*}
\end{prop}
\begin{proof}

The following is based on the calculations used in \cite[Lemma~3.2 and Proposition~3.3]{AW}
and \cite[Proposition~7.5]{AHRW}. Since $m_\tau$ is convex,
continuously differentiable and vanishes on $(-\infty,0]$,
\cite[Lemma~3.2]{AW} yields  $m_\tau(r)\leq r\,m_\tau'(r)$ for all $r\in\mathbb R$. Hence, pointwise in $D$,
\[
\left(
\frac1\tau m_\tau(y_\tau^*-\psi)|\bar{p}_\tau|
\right)^2
\leq
\left(
\frac1\tau m_\tau(y_\tau^*-\psi)(y_\tau^*-\psi)
\right)
\left(
\frac1\tau m_\tau'(y_\tau^*-\psi)|\bar{p}_\tau|^2
\right).
\]
Since $\zeta_\tau^*
=
\frac1\tau m_\tau(y_\tau^*-\psi)$ and
$\bar{\mu}_\tau
=
\frac1\tau m_\tau'(y_\tau^*-\psi)\bar{p}_\tau$,
the above reads, pointwise a.e. in $D$,
$(\zeta_\tau^* |\bar{p}_\tau|)^2 \leq \zeta_\tau^*(y_\tau^*-\psi)\bar{\mu}_\tau\bar{p}_\tau$. 
Taking square roots and integrating over $D$, the Cauchy--Schwarz inequality gives
\[
\langle\zeta_\tau^*,|\bar{p}_\tau|\rangle
\leq
\langle\zeta_\tau^*,y_\tau^*-\psi\rangle^{1/2}
\langle\bar{\mu}_\tau,\bar{p}_\tau\rangle^{1/2}.
\]
Now, testing the penalised adjoint equation \eqref{eq:ptau-old} with $\bar{p}_\tau$ gives a.s.,
\begin{equation}\label{eq:pTauPointwiseBound}
\|\bar{p}_\tau\|_V\leq C_a^{-1}\|J_y(y_\tau^*,\cdot)\|_{V^*},
\end{equation}
and using the non-negativity of $m_\tau'$,
$\langle\bar{\mu}_\tau,\bar{p}_\tau\rangle\geq0$. 
Using $\bar{\mu}_\tau = J_y(y_\tau^*)-A^*\bar{p}_\tau$ and the above two inequalities yields %
\[
0\leq
\langle\bar{\mu}_\tau,\bar{p}_\tau\rangle
\leq
C_a^{-1}\|J_y(y_\tau^*,\cdot)\|_{V^*}^2.
\]
Thus we have shown
\begin{equation}\label{eq:zetaAbsPPointwise}
0\leq
\langle\zeta_\tau^*,|\bar{p}_\tau|\rangle
\leq
C_a^{-1/2}
\langle\zeta_\tau^*,y_\tau^*-\psi\rangle^{1/2}
\|J_y(y_\tau^*,\cdot)\|_{V^*}
\qquad\text{a.s.}
\end{equation}
Now, 
      note that
\[
\begin{aligned}
&\langle\zeta_\tau^*,y_\tau^*-\psi\rangle
-\langle\zeta^*,y^*-\psi\rangle =
\langle\zeta_\tau^*-\zeta^*,y_\tau^*-\psi\rangle
+\langle\zeta^*,y_\tau^*-y^*\rangle,
\end{aligned}
\]
and both terms on the right-hand side converge to zero in $L^{q/2}(\Omega)$. Noticing that  $
\langle\zeta^*,y^*-\psi\rangle=0
$ and
$
\langle\zeta_\tau^*,y_\tau^*-\psi\rangle\geq0,
$
we obtain 
\begin{equation*} 
\langle\zeta_\tau^*,y_\tau^*-\psi\rangle^{1/2}
\to0
\qquad\text{in }L^q(\Omega).
\end{equation*}
Plugging this back into \eqref{eq:zetaAbsPPointwise} after taking the expectation,  and using boundedness of $J_y(y_\tau^*,\cdot)$ in $L^{\tilde p}(\Omega;V^*)$,
\[
0\leq\E[\langle\zeta_\tau^*,|\bar{p}_\tau|\rangle]
\leq C
\big\|\langle\zeta_\tau^*,y_\tau^*-\psi\rangle^{1/2}\big\|_{L^q(\Omega)}
\|J_y(y_\tau^*,\cdot)\|_{L^{\tilde p}(\Omega; V^*)}
\to0.
\]
Since $\zeta_\tau^*\to\zeta^*$ in $L^q(\Omega;V^*)$ and $(\bar{p}_\tau)$ is bounded in $L^{\tilde p}(\Omega;V)$, writing $\E[\langle\zeta^*,|\bar{p}_\tau|\rangle] = \E[\langle\zeta^*-\zeta_\tau^*,|\bar{p}_\tau|\rangle] + \E[\langle\zeta_\tau^*,|\bar{p}_\tau|\rangle]$ and using the above result we find
\begin{equation}
    \E[\langle\zeta^*,|\bar{p}_\tau|\rangle]\to0.\label{eq:new_0}
\end{equation}
Now, observe that the functional 
$v\mapsto\E[\langle\zeta^*,|v|\rangle]$
is weakly lower semicontinuous on $L^{\tilde p}(\Omega;V)$. Indeed it is continuous and as it is convex, also weakly lower semicontinuous. Thus, from $\bar{p}_\tau\weaklyto \bar{p}$, \eqref{eq:new_0} and using additionally the non-negativity of $\mathbb{E}[\langle \zeta^*, |\bar p| \rangle]$, we get
\[\E[\langle \zeta^*, |\bar{p}| \rangle] = 0.\]
Therefore \eqref{eq:pAbsOrthogonality}  holds since $\langle \zeta^*, |\bar{p}|\rangle$ is non-negative a.s. as a product of non-negative elements, so the previous equality is also zero on an a.s. level.

For the second claim, \eqref{eq:zetaAbsPPointwise} and H\"older's inequality yield
\[
\begin{aligned}
0\leq\E[\langle\zeta_\tau^*,|q_\tau^*|\rangle]
&\leq C\E\!\left[
\langle\zeta_\tau^*,y_\tau^*-\psi\rangle^{1/2}
|\pi_\tau^*|\|J_y(y_\tau^*,\cdot)\|_{V^*}
\right] \leq C
\big\|\langle\zeta_\tau^*,y_\tau^*-\psi\rangle^{1/2}\big\|_{L^q}
\|\pi_\tau^*\|_{L^{p'}}
\|J_y(y_\tau^*,\cdot)\|_{L^{\tilde p}}
\to0,
\end{aligned}
\]
where $1/q+1/p'+1/\tilde p=1$ from \eqref{eq:qprimeRelatedToPprimeEtc}. Since $\zeta_\tau^*\to\zeta^*$ in $L^q(\Omega;V^*)$ and $(q_\tau^*)$ is bounded in $L^{q'}(\Omega;V)$, the same argument as above gives \eqref{eq:qAbsOrthogonality}. The remaining positive- and negative-part identities follow from $|v|=v^++v^-$ and $\zeta^*\geq0$.
\end{proof}
Next, we shall strengthen \eqref{eq:eaEaQ} and \eqref{eq:eaEa} considerably by removing the $\mathcal{E}$-almost part (which arose from an application of Egorov's theorem) by utilising recently-developed improved analysis \cite{AHRW}. Recall the coincidence set $\mathcal{A}^*(\omega)$ defined in \eqref{eq:coincidence_set}.
\begin{prop}\label{prop:exactInactiveAnnihilation}
We have for a.s. $\omega \in \Omega$,
\begin{align}
\langle \lambda^*(\omega),v\rangle&=0
\qquad\forall v\in V:\ v=0\quad\text{q.e. on }\mathcal A^*(\omega),\label{eq:qResidualExactSupport}\\
\langle \bar{\mu}(\omega),v\rangle&=0
\qquad\forall v\in V:\ v=0\quad\text{q.e. on }\mathcal A^*(\omega).\label{eq:pResidualExactSupport}
\end{align}
\end{prop}
\begin{proof}
We follow the approach in \cite[Section 7]{AHRW}. 
Set
\[
s_\tau:=(y_\tau^*-\psi)^-,
\qquad
s:=(y^*-\psi)^-=\psi-y^*.
\]
The continuity of the negative-part mapping and a Dominated Convergence Theorem argument give $s_\tau\to s$ in $L^q(\Omega;V)$ and therefore also in $L^{\tilde p'}(\Omega;V)$.

Take an arbitrary $v\in V$ with $0\leq v\leq s_\tau(\omega)$. We note that on $\{y_\tau^*-\psi\leq0\}$ the term $m_\tau'(y_\tau^*-\psi)$ vanishes, while on $\{y_\tau^*-\psi>0\}$ the inequality $v\leq s_\tau$ forces $v=0$. This implies that for a.s. $\omega$ we have
\begin{equation}\label{eq:penalisedIntervalAnnihilation}
\langle {\lambda}_\tau^*(\omega),v\rangle = 
\pi_\tau^*(\omega)\langle\bar{\mu}_\tau(\omega),v\rangle=0.
\end{equation}
If $v\in L^{\tilde p'}(\Omega;V)$ satisfies $0\leq v\leq s$, define, as in \cite[Lemma 7.6]{AHRW},
\[
v_\tau:=\min(v,s_\tau).%
\]
Then $v_\tau\to v$ in $L^{\tilde p'}(\Omega;V)$, $0\leq v_\tau\leq s_\tau$, and hence by \eqref{eq:penalisedIntervalAnnihilation} and the weak convergence of $\bar{\mu}_\tau$ in \eqref{eq:lambdaWeakLimits}, we get
\[
0=\E[\langle\bar{\mu}_\tau,v_\tau\rangle]
\to\E[\langle \bar{\mu},v\rangle],
\]
i.e.,
\begin{equation}\label{eq:pResidualIntervalIntegrated}
\E[\langle \bar{\mu},v\rangle]=0
\qquad\forall v\in L^{\tilde p'}(\Omega;V):\ 0\leq v\leq s.
\end{equation}
Now we want to pass from this into an a.s. statement. We do so with a standard localisation argument. For this purpose, let $\{d_j\}_{j\in\mathbb N}$ be dense in $V$ and set
\[
v_j(\omega):=\min(d_j^+,s(\omega)).
\]
Note that $v_j \in L^{\tilde p'}(\Omega;V)$ and $0\leq v_j\leq s$.
Hence, for every $E\in\mathcal F$, the function $\chi_{E}v_j$ is admissible in
\eqref{eq:pResidualIntervalIntegrated}, and therefore
\[
\int_E \langle \bar{\mu}(\omega),v_j(\omega)\rangle\,\mathrm{d}\mathbb P(\omega)=0.
\]
Since this holds for every $E\in\mathcal F$, we must have
\[
\langle \bar{\mu}(\omega),v_j(\omega)\rangle=0
\qquad\text{for a.e. }\omega.
\]
Taking the union over $j\in\mathbb N$ of the exceptional sets associated to the above displayed equality, there is a
single null set $N$ such that the preceding identity holds for every $j$
and every $\omega\notin N$. Now we argue as in \cite[Proposition 7.8]{AHRW}. Fix $\omega\notin N$ and let $v\in V$ satisfy $0\leq v\leq s(\omega)$.
Choose a subsequence $(d_{j_k})$ with $d_{j_k}\to v$ in $V$. Then
\[
v_{j_k}(\omega)
=\min(d_{j_k}^+,s(\omega))
\to
\min(v,s(\omega))
=v
\qquad\text{in }V.
\]
It follows that
\[
\langle \bar{\mu}(\omega),v\rangle=0.
\]
Thus, outside the null set $N$,
\begin{equation}
\langle \bar{\mu}(\omega),v\rangle=0
\qquad
\forall v\in V:\ 0\leq v\leq s(\omega).
\label{eq:new_2}    
\end{equation}
From here, to reach the quasi-everywhere representation in \eqref{eq:pResidualExactSupport} is standard, but we give the details for completeness, arguing exactly as in the proof of \cite[Proposition 7.8]{AHRW}. Note that we have $\mathcal{A}^*(\omega) = \{ s(\omega) = 0 \}.$ Take $v \in V$ with $v \geq 0$ and $v = 0$ q.e. on $\mathcal{A}^*(\omega)$. Utilising the following expression for the tangent cone of $V_+$ (see \cite[Theorem 6.57]{MR1756264} for the case $V=H^1_0(D)$ or \cite[Lemme 3.2]{MR0423155} in the general Dirichlet space setting):
\[\mathcal{T}_{V_+}(s) = \{ \varphi \in V : \varphi \geq 0 \text{ q.e. on $\mathcal{A}^*(\omega)$} \},\]
it follows that $-v \in \mathcal{T}_{V^+}(s(\omega))$ and hence, there exists a sequence $\{v_n\}$ with $v_n \to v$ in $V$ and $v_n \le t_n s(\omega)$ for some $t_n > 0$. Thus,
    \[0 \le \max( 0, v_n / t_n ) \le s(\omega)\]
and we can apply \eqref{eq:new_2} and get
    \[ \langle \bar{\mu}(\omega), \max(0, v_n / t_n) \rangle = 0.\]
Multiplying by $t_n$ and passing to the limit $n \to \infty$ gives
  \[  \langle \bar{\mu}(\omega), v \rangle = 0   \text{ for all $v \ge 0$, $v = 0$ q.e. on $\mathcal A^*(\omega)$.}
\]
Using the decomposition $v=v^+ - v^-$, we obtain the second claim.

Repeating the same argument with $\lambda^*_\tau$, using \eqref{eq:weightedLambdaWeakLimit} and \eqref{eq:penalisedIntervalAnnihilation} and $s_\tau\to s$ in $L^q(\Omega;V)$, we retrieve also the first claim.
\end{proof}

Now we focus on obtaining $\langle \bar \mu, \bar p \rangle \geq 0$. It is useful to define $a(\omega; \cdot, \cdot)\colon V \times V \to \mathbb{R}$ as the bilinear form generated by $A(\omega)$, i.e.,
\[a(\omega; y, z) := \langle A(\omega)y, z \rangle.\]
From the penalised adjoint equation \eqref{eq:ptau-old} we obtain
\begin{equation}\label{eq:quadraticEnergyP}
\langle J_y(y_{\tau}^*,\omega),\bar{p}_{\tau}(\omega)\rangle-a(\omega; \bar{p}_{\tau}(\omega),\bar{p}_{\tau}(\omega))=\langle\bar \mu_{\tau}(\omega),\bar{p}_{\tau}(\omega)\rangle\ge0 .
\end{equation}
As we have only weak convergences of $\bar \mu_{\tau}$ and $\bar{p}_\tau$, we cannot directly pass to the limit (neither in the inequality above nor the one after taking expectations). We will argue using a Minty approach instead. The next proposition provides the tool to do so.

\begin{prop}\label{lem:technical_tool}
Assume $\theta_\tau \weaklyto \theta$ in $L^1(\Omega)$ where $\theta_\tau \in L^1(\Omega)$ is non-negative and suppose that $\theta_\tau\bar{p}_{\tau}  \weaklyto \eta$ to some $\eta$ in $L^1(\Omega; V)$. Then we have
\begin{equation}\label{eq:limitEnergy}
\theta(\omega)\langle J_y(y^*(\omega),\omega),\eta(\omega)\rangle-a(\omega;\eta(\omega),\eta(\omega))\ge0 \quad\text{a.s. in $\Omega$}.
\end{equation}
\end{prop}
\begin{proof}
The idea is to derive an affine inequality which allows us to use the weak convergences that we do have in combination with Mazur's lemma to pass to the limit, and then recover the desired statement. Let $\tau_k\downarrow0$ denote the subsequence
along which the convergences in the assumptions hold and set
\[
\eta_{\tau_k}:=\theta_{\tau_k}\bar p_{\tau_k}.
\]
\noindent\emph{Step 1: an affine inequality.}
For any $z\in V$, coercivity of $a(\omega;\cdot, \cdot)$ gives %
\[a(\omega;\bar{p}_{\tau_k}(\omega),\bar{p}_{\tau_k}(\omega))\ge a(\omega;\bar{p}_{\tau_k}(\omega),z)+a(\omega;z,\bar{p}_{\tau_k}(\omega))-a(\omega;z,z).\]
Together with
\eqref{eq:quadraticEnergyP},
\[
\langle J_y(y_{\tau_k}^*(\omega),\omega),\bar{p}_{\tau_k}(\omega)\rangle\ge a(\omega;\bar{p}_{\tau_k}(\omega),z)+a(\omega;z,\bar{p}_{\tau_k}(\omega))-a(\omega;z,z).
\]
Multiplying by $\theta_{\tau_k}\ge0$ yields
\begin{equation}\label{eq:affineIneq}
\langle J_y(y_{\tau_k}^*(\omega),\omega),\eta_{\tau_k}(\omega)\rangle\;\ge\;a(\omega;\eta_{\tau_k}(\omega),z)+a(\omega;z,\eta_{\tau_k}(\omega))-\theta_{\tau_k}\,a(\omega;z,z)
\qquad\text{for every }z\in V .
\end{equation}

\medskip
\noindent\emph{Step 2: passage to the limit.}
By the premise we have
$(\theta_{\tau_k},\eta_{\tau_k})\rightharpoonup(\theta,\eta)$ in $L^1(\Omega)\times L^1(\Omega;V)$. %
Applying Mazur's lemma to this pair gives for every $N$ the existence of a number $K_N \geq N$ and coefficients
$\alpha_{N,k} \ge 0$ for $k=N,\dots,K_N$ with $\sum_{k=N}^{K_N}\alpha_{N,k}=1$ such that
\[
\bar\theta_N:=\sum_{k=N}^{K_N}\alpha_{N,k}\theta_{\tau_k},\qquad
\bar \eta_N:=\sum_{k=N}^{K_N}\alpha_{N,k}\eta_{\tau_k},\]
and
\[(\bar\theta_N,\bar \eta_N)\to(\theta,\eta)\ \text{ strongly in }L^1(\Omega) \times L^1(\Omega;V).
\]
After passing to a further subsequence, %
$\bar\theta_N(\omega)\to\theta(\omega)$ and $\bar \eta_N(\omega)\to \eta(\omega)$ in $V$ for almost every $\omega$.

In \eqref{eq:affineIneq}, writing $J_y(y_{\tau_k}^*,\omega) = J_y(y^*, \omega) + J_y(y_{\tau_k}^*,\omega)-J_y(y^*, \omega)$, multiplying by $\alpha_{N,k}$ and summing yields
\begin{align}\label{eq:affineMazur}
\langle J_y(y^*(\omega),\omega),\bar \eta_N(\omega)\rangle+E_N(\omega)
\;\ge\;a(\omega;\bar \eta_N(\omega),z)+a(\omega;z,\bar \eta_N(\omega))-\bar\theta_N\,a(\omega;z,z),
\end{align}
where the error term is
\[E_N(\omega):=\sum_{k=N}^{K_N}\alpha_{N,k}\langle J_y(y_{\tau_k}^*(\omega),\omega)-J_y(y^*(\omega),\omega),\eta_{\tau_k}(\omega)\rangle.\]
Since $\eta_{\tau_k}=\theta_{\tau_k}\bar{p}_{\tau_k}$ with $\theta_{\tau_k}\ge0$, the pointwise bound \eqref{eq:pTauPointwiseBound} on $\bar{p}_{\tau_k}$
gives
\[
|E_N(\omega)|\le C_a^{-1}\sup_{k\ge N}\|J_y(y_{\tau_k}^*(\omega),\omega)-J_y(y^*(\omega), \omega)\|_{V^*}\,\sup_{k\ge N}\|J_y(y_{\tau_k}^*(\omega),\omega)\|_{V^*}\;\bar\theta_N(\omega) .
\]
In the limit $N \to \infty$, the first factor tends to $0$ (as it is the convergence of the tail of a convergent sequence), the second is bounded (as the sequence is convergent), and as
$\bar\theta_N(\omega)\to\theta(\omega)$ the third term is also bounded, thus $E_N(\omega)\to0$. Letting $N\to\infty$ in \eqref{eq:affineMazur} gives therefore
\begin{equation}\label{eq:affineLimit}
\langle J_y(y^*,\omega),\eta(\omega)\rangle\;\ge\;a(\omega;\eta(\omega),z)+a(\omega;z,\eta(\omega))-\theta(\omega)a(\omega;z,z)
\qquad \text{ for a.e. } \omega, \; \forall z\in V.
\end{equation}

\medskip
\noindent\emph{Step 3: choice of $z$.} Now we select $z$ carefully in order to obtain \eqref{eq:limitEnergy}.

\emph{Case $\theta(\omega)>0$.} Take $z=\eta(\omega)/\theta(\omega)$ in \eqref{eq:affineLimit}. Then
$a(\omega;\eta,z)+a(\omega;z,\eta)=2a(\omega;\eta,\eta)/\theta$ and $\theta\,a(\omega;z,z)=a(\omega;\eta,\eta)/\theta$, so
\[
\langle J_y(y^*(\omega),\omega),\eta(\omega)\rangle\ge\frac{a(\omega;\eta(\omega),\eta(\omega))}{\theta(\omega)}.
\]
Multiplying by $\theta$ gives \eqref{eq:limitEnergy}.

\emph{Case $\theta(\omega)=0$.} Then \eqref{eq:affineLimit} reads
$\langle J_y(y^*,\omega),\eta(\omega)\rangle\ge a(\omega;\eta(\omega),z)+a(\omega;z,\eta(\omega))$ for all $z$. Choosing $z=M\eta$ for some arbitrary $M>0$ gives
$\langle J_y(y^*(\omega),\omega),\eta(\omega)\rangle\ge 2M\,a(\omega;\eta(\omega),\eta(\omega))$; dividing by $M$ and letting $M\to\infty$ forces $a(\omega;\eta,\eta)\le0$, %
so \eqref{eq:limitEnergy} holds (trivially).
\end{proof}
Using this result, we can first strengthen \eqref{eq:TT} to an a.s. statement.
\begin{prop}\label{prop:lambda_p_nonneg}
    We have
    \begin{equation}\label{eq:pstarEnergy}
\langle \bar{\mu},\bar{p}\rangle\geq0,
\qquad\text{a.s.}
\end{equation}
\end{prop}
\begin{proof}
    We apply \cref{lem:technical_tool} with the choice $\theta_\tau \equiv 1$ in which case note that $\eta = \bar{p}$ due to \cref{lem:bounded-dual-variables}. Recalling the definition of $\bar{\mu}$, we get the result immediately.
\end{proof}

We now state the following intermediary stationarity result, which subsumes and improves the systems \eqref{eq:eaWeakStationarity} and \eqref{eq:ecsSytem}. In addition, we do it under weaker circumstances (in \cite[Assumption 3.8]{AGHS}, it was additionally assumed that, for every
$v\in U_{ad}$, the map
$u\mapsto\varrho'(u)(u-v)$ is weakly lower semicontinuous; as we saw in \cref{prop:new_1}, this is not necessary).
\begin{prop}[Enhanced weak stationarity]\label{thmenhancedWeakStationarity}
Let \cref{ass:newCombinedOnAandB,ass:onJ,ass:onvarrho} hold, and let
$u^*\in U_{ad}$ be a local minimiser of \eqref{eq:OCProblem},
with $y^*=S(u^*)$. Then there exist
\[
\zeta^*\in L^q(\Omega;V^*),\qquad
q^*\in L^{q'}(\Omega;V),\qquad
\lambda^*\in L^{q'}(\Omega;V^*),\qquad
\pi^*\in L^{p'}(\Omega),
\]
such that
\begin{subequations}\label{eq:enhancedWeakSystem}
\begin{align}
A(\omega)y^*(\omega)-f(\omega)-B(\omega)u^*
+\zeta^*(\omega)
&=0 &&\text{a.s.},
\\
\zeta^*(\omega)\ge0,\quad
y^*(\omega)\le\psi(\omega),\quad
\langle\zeta^*(\omega),y^*(\omega)-\psi(\omega)\rangle
&=0 &&\text{a.s.},
\\
A(\omega)^*q^*(\omega)+\lambda^*(\omega)
-\pi^*(\omega)J_y(y^*(\omega),\omega)&=0
&&\text{a.s.},
\\
\E[\langle B^*q^*,u^*-v\rangle_{U^*,U}]
+\varrho'(u^*)(u^*-v)
&\le0 &&\forall v\in U_{ad},
\\
\pi^*&\in\partial\cR[\cJ(y^*)],
\\
\langle\zeta^*(\omega),(q^*(\omega))^+\rangle
=\langle\zeta^*(\omega),(q^*(\omega))^-\rangle
&=0 &&\text{a.s.},
\\
\langle\lambda^*(\omega),v\rangle
&=0
&& 
\forall v\in V:\ v=0\ \text{q.e. on }\mathcal A^*(\omega), 
\text{ a.s.},%
\end{align}
\end{subequations}
In addition, there exist
\[
\bar{p}\in L^{\tilde p}(\Omega;V),\qquad
\bar{\mu}\in L^{\tilde p}(\Omega;V^*),\qquad
\]
such that
\begin{subequations} 
\begin{align}
A(\omega)^*\bar{p}(\omega)+\bar{\mu}(\omega)
-J_y(y^*(\omega),\omega)&=0
&&\text{a.s.},\label{eq:system_barp_eqn}
\\
\langle\zeta^*(\omega),(\bar{p}(\omega))^+\rangle
=\langle\zeta^*(\omega),(\bar{p}(\omega))^-\rangle
&=0 &&\text{a.s.},
\\
\langle\bar{\mu}(\omega),v\rangle
&=0
&& 
\forall v\in V:\ v=0\ \text{q.e. on }\mathcal A^*(\omega), 
\text{ a.s.},\\
\langle\bar{\mu}(\omega),\bar{p}(\omega)\rangle
&\ge0 &&\text{a.s.}\label{eq:system_barp_barmu_sign}
\end{align}
\end{subequations}
\end{prop}
\begin{proof}
    Putting together \eqref{eq:zetaStrongNew}, \eqref{eq:lambdaWeakLimits}, \eqref{eq:weightedLambdaWeakLimit} that define $\zeta^*, \bar\mu$ and $\lambda^*$, \cref{prop:new_1} that relates $q^*$ and $u^*$, as well as the new conditions derived in 
\cref{prop:separateContactOrthogonality,prop:exactInactiveAnnihilation,prop:lambda_p_nonneg}, we obtain the result for any minimiser $u^*$ arising as a limit of penalised minimisers as in \cref{lem:boundednessofyandutaus}. To upgrade from this to  the result for arbitrary local minimisers as claimed, we use a classical localisation trick, see e.g. \cite[\S 5]{AHROCQVI} for closely-related arguments. For a local
minimiser $u^*$, choose $\varepsilon>0$ such that $u^*$
minimises the original objective over
$U_{ad}\cap\overline B_\varepsilon(u^*)$.
We then consider the localised penalised problems
\[
\min_{u\in U_{ad}\cap\overline B_\varepsilon(u^*)}
\cR[\cJ(T_\tau(u))]+\varrho(u)
+\frac12\|u-u^*\|_U^2.
\]
Firstly, note that since $U_{ad}\cap\overline B_\varepsilon(u^*)$ is bounded, the localised problems admit minimisers and all needed results in the previous section hold for them without any assumptions regarding the uniform coercivity of $\{\mathcal R\circ\mathcal J\circ T_\tau+\varrho\}_{\tau>0}$ (see the discussion below \eqref{eq:penalisedOCProblem}).

Now, by comparison with $u^*$, the convergence properties
of $T_\tau$, and weak lower semicontinuity, the
corresponding minimisers $u_\tau^*$ converge strongly
to $u^*$ in $U$. Thus, for sufficiently small
$\tau$, the ball constraint is inactive. The penalised
optimality condition therefore holds for all
$v\in U_{ad}$, with the additional term
\[
(u_\tau^*-u^*,u_\tau^*-v)_U,
\]
which vanishes as $\tau\to0$. All remaining estimates
and limit arguments carry over unchanged, yielding
the asserted stationarity system at $u^*$.
\end{proof}
Comparing this result to our main result \cref{thmmain}, we see that we obtain all of \eqref{eq:factorisedSystem} except the sign condition \eqref{eq:lambdaqgeq0Nonsmooth} which is not shown, and we obtain a system similar to \eqref{eq:system_p_star} but not the VI \eqref{eq:viHat} relating the control to the adjoint state, nor any relationship between $\bar p$, $\pi^*$ and $q^*$. Although \cref{thmenhancedWeakStationarity} provides
an unweighted adjoint $\bar p$ satisfying the a.s. 
C-stationarity conditions, the control optimality
condition involves the weighted adjoint $q^*$.
In general, the weak limits do not satisfy
$q^*=\pi^*\bar p$, and the result does not guarantee
the existence of an unweighted adjoint $p^*$ such that
$q^*=\pi^*p^*$ and the scenario-wise C-stationarity
conditions hold simultaneously. Thus, the system
does not yet constitute full C-stationarity.
We establish the required compatibility under
additional assumptions on $\cR$ in the following two sections.

\begin{remark}
    From now on, $(u_\tau^*,y_\tau^*,\pi_\tau^*,\bar p_\tau,q_\tau^*)$ and their limits refer to the localised penalised problems. Since the arguments of this section and the remaining sections use only the penalised state and adjoint equations and $\pi_\tau^*\in\partial\mathcal R[\mathcal J(y_\tau^*)]$, and not the control inequality, all subsequent results apply with $u^*$ the given local minimiser.
\end{remark}

\section{C-stationarity for Fréchet differentiable risk measures}\label{sec:C_stat_Frechet}
In this section, we will prove the main result \cref{thmmain} under Fréchet differentiability of $\cR$ at $\cJ(y^*)$. We begin with the following, which improves \cite[Proposition 4.11]{AGHS} as it needs smoothness only at a single point.
\begin{prop}\label{prop:identification}
Suppose that $\cR$ is Fr\'echet differentiable at
$\cJ(y^*)$. Then
\[
\pi_\tau^*\to\pi^*
\text{ in }L^{p'}(\Omega)
\qquad\text{and}\qquad q^*=\pi^*\bar p.
\]
\end{prop}

\begin{proof}
Set $\cJ_\tau:=\cJ(y_\tau^*)$ and
$\cJ^*:=\cJ(y^*)$. We already know that
$\cJ_\tau\to \cJ^*$ in $L^p(\Omega)$ and
$\pi_\tau^*\in\partial\cR[\cJ_\tau]$. Since $\cR$ is convex and Fr\'echet differentiable at $\cJ^*$, its subdifferential at this point is the singleton $\{\cR'[\cJ^*]\}$. Thus $\pi^*=\cR'[\cJ^*]$.

We first prove that $\pi_\tau^*\to\pi^*$
strongly in $L^{p'}(\Omega)$. 
Define 
\[
h_\tau:=\cJ_\tau-\cJ^*,
\qquad
r(h):=\cR[\cJ^*+h]-\cR[\cJ^*]-\E[\pi^*h].
\]
By Fr\'echet differentiability, for any $\varepsilon>0$, there exists
$\delta>0$ such that
\[
|r(h)|\leq\varepsilon\|h\|_{L^p(\Omega)}
\qquad
\text{whenever }\|h\|_{L^p(\Omega)}\leq\delta.
\]
Fix $0<t<\delta/2$ and take $\tau$ sufficiently
small that $\|h_\tau\|_{L^p(\Omega)}\leq t/2$.
Then, for every $v\in L^p(\Omega)$ with
$\|v\|_{L^p(\Omega)}=1$, the above remainder estimate
applies to $h_\tau$ and $h_\tau\pm tv$.

Using $\pi_\tau^*\in\partial\cR[\cJ_\tau]$,
we obtain
\begin{align*}
\E[(\pi_\tau^*-\pi^*)v]
&\leq
\frac{\cR[\cJ_\tau+tv]-\cR[\cJ_\tau]}{t}
-\E[\pi^*v]\\
&=\frac{r(h_\tau+tv)-r(h_\tau)}{t}\\
&\leq
\frac{\varepsilon}{t}
\bigl(\|h_\tau+tv\|_{L^p(\Omega)}
+\|h_\tau\|_{L^p(\Omega)}\bigr)\\
&\leq
\varepsilon\left(
1+\frac{2\|h_\tau\|_{L^p(\Omega)}}{t}
\right)\\
&\leq2\varepsilon.
\end{align*}
Replacing $v$ by $-v$ gives the same bound
for $-\E[(\pi_\tau^*-\pi^*)v]$. Hence
\[
\left|\E[(\pi_\tau^*-\pi^*)v]\right|
\leq2\varepsilon
\qquad
\forall v\in L^p(\Omega):
\|v\|_{L^p(\Omega)}=1.
\]
Taking the supremum over these $v$ yields
\[
\|\pi_\tau^*-\pi^*\|_{L^{p'}(\Omega)}
\leq2\varepsilon.
\]
Since $\varepsilon>0$ was arbitrary, we find $
\pi_\tau^*\to\pi^*$. The rest of the proof is basically identical to that of \cite[Proposition 4.11]{AGHS}. We give it for completeness. Recall that we have
\[
\bar{p}_\tau\rightharpoonup\bar p
\quad\text{in }L^{\tilde p}(\Omega;V),
\qquad
\pi_\tau^*\bar{p}_\tau\rightharpoonup q^*
\quad\text{in }L^{q'}(\Omega;V).
\]
For any $\varphi\in L^q(\Omega;V^*)$,
H\"older's inequality and $\frac1{\tilde p'}=\frac1{p'}+\frac1q$ from \eqref{eq:qprimeRelatedToPprimeEtc}
give
\[
\|(\pi_\tau^*-\pi^*)\varphi\|_
{L^{\tilde p'}(\Omega;V^*)}
\leq
\|\pi_\tau^*-\pi^*\|_{L^{p'}(\Omega)}
\|\varphi\|_{L^q(\Omega;V^*)}
\to0.
\]
Hence,
\[
\begin{aligned}
\E[\langle q^*,\varphi\rangle]
&=\lim_{\tau\downarrow0}
  \E[\langle\pi_\tau^*\bar{p}_\tau,\varphi\rangle]
=\lim_{\tau\downarrow0}
  \E[\langle \bar{p}_\tau,\pi_\tau^*\varphi\rangle]
=\E[\langle\bar p,\pi^*\varphi\rangle]
=\E[\langle\pi^*\bar p,\varphi\rangle].
\end{aligned}
\]
Since $\varphi$ was arbitrary, we conclude
that $q^*=\pi^*\bar p$.
\end{proof}
Now we come to proving one branch of \cref{thmmain}.
\begin{proof}[\bf Proof of \cref{thmmain} under Fréchet differentiability of $\cR$] By \cref{thmenhancedWeakStationarity}, there exist
$\zeta^*$, $\pi^*$, $\bar p$, $q^*$,
$\bar\mu$ and $\lambda^*$ satisfying the enhanced weak stationarity system \eqref{eq:enhancedWeakSystem}.  As $\cR$ is Fr\'echet differentiable at
$\cJ(y^*)$, by \cref{prop:identification} its subdifferential at this point
is the singleton
\[
\partial\cR[\cJ(y^*)]
=\{\cR'[\cJ(y^*)]\},
\]
thus $\pi^*=\cR'[\cJ(y^*)]$, and we have the identification
\[
q^*=\pi^*\bar p
\qquad\text{a.s.}
\]
We may therefore set
\[
p^*:=\bar p,
\qquad
\mu^*:=\bar\mu .
\]
Hence $q^* = \pi^*p^*$.  Moreover, by the definitions of $\lambda^*$ and $\bar\mu$,
\[
\lambda^*
=\pi^*J_y(y^*,\cdot)-A^*q^*
=\pi^*(J_y(y^*,\cdot)-A^*p^*)
=\pi^*\mu^*
\qquad\text{a.s.}
\]
The statement \eqref{eq:lambdaqgeq0Nonsmooth} follows due to this identification and \eqref{eq:system_barp_barmu_sign}:
\begin{align*}
    \langle \lambda^*(\omega), q^*(\omega) \rangle
&=\langle \pi^*(\omega)J_y(y^*(\omega),\omega) -  A(\omega)^*q^*(\omega), q^*(\omega) \rangle\\
&=(\pi^*(\omega))^2\langle J_y(y^*(\omega),\omega) - A(\omega)^*p^*(\omega), p^*(\omega) \rangle\\
&= (\pi^*(\omega))^2\langle \mu^*(\omega), p^*(\omega) \rangle \geq 0.
\end{align*}

\end{proof}

\section{C-stationarity for monotone risk measures}\label{sec:factorisation}
We now show that C-stationarity holds also for (potentially nonsmooth) monotone risk measures, i.e., we prove \cref{thmmain} under monotonicity of $\cR$. First, we need an important consequence of the monotonicity.
\begin{lemma}[Non-negativity of risk identifiers]\label{lem:piNonnegative}
Let $\cR$ be monotone. For $X \in L^p(\Omega)$, every $\pi\in\partial\cR[X]$ satisfies $\pi\geq0$ a.s. In particular,
\[
\pi_\tau^*\geq0,\qquad \pi^*\geq0\qquad\text{a.s.}
\]
\end{lemma}
\begin{proof}
Let $Z\in L^p(\Omega)$ with $Z\geq0$. Using the monotonicity of $\cR$ and by definition of the subdifferential,%
\[
0\leq \cR[X]-\cR[X-Z]\leq \E[\pi Z].
\]
Since this holds for every non-negative $Z$, we find $\pi\geq0$ a.s. 
\end{proof}
Now using the monotonicity of the risk measure, we can improve the inconvenient limsup condition \eqref{eq:limsupCond} to an a.s. condition stated for the limit. Unlike in the proof of the theorem in \cref{sec:C_stat_Frechet} which obtained the condition \eqref{eq:weightedResidualSign} below from the equivalent condition for $p^*$, we do it in reverse, since we do not yet have the identification of $q^*$ and $\pi^*p^*$, which was needed for that.
\begin{prop}\label{prop:qstar_complementarity}
Let $\cR$ be monotone. Then
\begin{align}
\langle \lambda^*,q^*\rangle &\geq 0
\quad\text{a.s.}\label{eq:weightedResidualSign},\\
q^* &= 0\quad\text{on }\{\pi^*=0\}.\label{eq:qstarEnergy}
\end{align}
\end{prop}
\begin{proof}
    In this case we apply \cref{lem:technical_tool} with $\theta_\tau := \pi_{\tau}^*$ so that $\theta = \pi^*$ and $\eta = q^*$. Thanks to the monotonicity of $\cR$,  $\theta_\tau \geq 0$, so the result is applicable and we obtain
\[\pi^*(\omega)\langle J_y(y^*(\omega),\omega),q^*(\omega)\rangle-a(\omega;q^*(\omega),q^*(\omega))\ge0 \quad\text{a.s. in $\Omega$}.\]
Recalling the definition of $\lambda^*$ we obtain the first statement. If $\pi^*(\omega)=\theta(\omega)$ vanishes, then the above displayed inequality and the coercivity of $a$ imply that  $q^*(\omega)=0$. 
\end{proof}
This gives a complementarity-type relation
between the risk identifier $\pi^*$ and the
risk-weighted adjoint $q^*$: the latter vanishes
on scenarios to which the former assigns zero weight. 
This observation suggests recovering an unweighted
adjoint from $q^*$. As the weak limit
of $\pi_\tau^*\bar p_\tau$ need not coincide with
$\pi^*\bar p$, the auxiliary adjoint $\bar p$
cannot in general be used directly in the
risk-weighted control optimality condition.

To obtain a C-stationarity system in terms of an
unweighted adjoint $p^*$, we seek a factorisation
$q^*=\pi^*p^*$. On $\{\pi^* > 0\}$, this necessarily
gives $p^*=q^*/\pi^*$. By
\cref{prop:qstar_complementarity}  we have indeed $q^*=0$ on
$\{\pi^*=0\}$, so the factorisation remains valid
regardless of the choice of $p^*$ on this set.
We define
\begin{equation}
    \label{eq:defn_of_pstar_new}
p^*:=
\begin{cases}
q^*/\pi^* &\text{if } \pi^*>0,\\
\bar p &\text{if }\pi^*=0.
\end{cases}
\end{equation}
The choice of $\bar p$ on $\{\pi^*=0\}$ allows us
to inherit the adjoint and multiplier properties already established for $\bar p$, which, combined
with those derived from $q^*$ on $\{\pi^*>0\}$,
will yield the full C-stationarity system.

We finally come to proving the remaining branch of  \cref{thmmain} and consequently finishing the proof.

\begin{proof}[\bf Proof of \cref{thmmain} under monotonicity of $\cR$]
Again by \cref{thmenhancedWeakStationarity}, there exist
$\zeta^*$, $\pi^*$, $\bar p$, $q^*$,
$\bar\mu$ and $\lambda^*$ satisfying the enhanced weak stationarity system \eqref{eq:enhancedWeakSystem}.  
As explained above, by \cref{prop:qstar_complementarity}, $q^*=0$ on $\{\pi^*=0\}$. Thus define $p^*$ as in \eqref{eq:defn_of_pstar_new}, i.e.,
\[
p^*:=
\begin{cases}
q^*/\pi^*,&\pi^*>0,\\
\bar{p},&\pi^*=0,
\end{cases}
\]
which by construction satisfies $\pi^*p^* = q^*$. 

Together with the bounds %
\begin{align*}
\|\bar{p}(\omega)\|_V &\leq C_a^{-1}\|J_y(y^*(\omega),\omega)\|_{V^*},
\\
\|q^*(\omega)\|_V&\leq C_a^{-1}\pi^*(\omega)\|J_y(y^*(\omega),\omega)\|_{V^*}
\quad\text{on }\{\pi^*>0\}    
\end{align*}
(these follow from the equations for $\bar{p}$ and $q^*$ in \eqref{eq:system_barp_eqn} and \eqref{eq:system_qstar_eqn} and the sign conditions \eqref{eq:system_barp_barmu_sign} and \eqref{eq:weightedResidualSign} respectively)
we deduce
\[
\|p^*(\omega)\|_V
\leq C_a^{-1}\|J_y(y^*(\omega),\omega)\|_{V^*}
\qquad\text{a.s.},
\]
giving $p^*\in L^{\tilde p}(\Omega;V)$. Set
\[
\mu^*:=J_y(y^*,\cdot)-A^*p^*.
\]
Then $\mu^*\in L^{\tilde p}(\Omega;V^*)$ and \eqref{eq:adjHat} holds. 
Furthermore, as by construction $q^*=\pi^*p^*$,
\[
\lambda^*
=\pi^*J_y(y^*,\cdot)-A^*q^*
=\pi^*(J_y(y^*,\cdot)-A^*p^*)
=\pi^*\mu^*
\qquad\text{a.s.}
\]
In addition, \eqref{eq:viHat} is the control inequality \eqref{eq:system_q_star_control_ineq} rewritten in terms of the single adjoint $p^*$.

On $\{\pi^*>0\}$,
\[
\mu^*
=\frac{\pi^*J_y(y^*,\cdot)-A^*q^*}{\pi^*},
\]
and therefore \eqref{eq:weightedResidualSign} gives
\[
\langle\mu^*,p^*\rangle
=\frac{\langle\pi^*J_y(y^*,\cdot)-A^*q^*,q^*\rangle}{(\pi^*)^2} = 
\frac{\langle\lambda^*,q^*\rangle}{(\pi^*)^2} \geq0\qquad\text{on $\{\pi^*>0\}$}.
\]
On $\{\pi^*=0\}$, $p^*=\bar{p}$ and
$\mu^*=J_y(y^*,\cdot)-A^*\bar{p} = \bar{\mu}$, so \eqref{eq:pstarEnergy} gives the same conclusion for $\mu^*$ and $p^*$, i.e., we have \eqref{eq:signHat}.

Likewise, on $\{\pi^*>0\}$ we have $(q^*)^\pm=\pi^*(p^*)^\pm$, whereas on $\{\pi^*=0\}$ we have $p^*=\bar{p}$.  \cref{prop:separateContactOrthogonality} therefore gives \eqref{eq:zetaHatPM}.

Finally, on $\{\pi^*>0\}$ the multiplier $\mu^*$ is the weighted residual in \eqref{eq:qResidualExactSupport} divided by $\pi^*$, whereas on $\{\pi^*=0\}$ it is the residual in \eqref{eq:pResidualExactSupport}. \cref{prop:exactInactiveAnnihilation} thus gives \eqref{eq:exactSupportHat}. %
\end{proof}

 \section{Conclusion}\label{sec:conclusion}
 In this paper, we have derived an optimality system for the risk-averse stochastic optimal control problem \eqref{eq:OCProblem} by utilising some recent developments in the control of obstacle problems and a novel approach related to designing a factorisation involving the adjoints, filling a gap in the literature. We leave other interesting aspects such as developing applications of the theory and numerical methods for solving such nonsmooth risk-averse problems for future work.

\bibliographystyle{abbrv}
\bibliography{references}
\end{document}